\documentclass[MSNbibl,number,citesort,dvips]{arxbj}
\usepackage{mathbh}

\aid{0}
\volume{20}
\issue{4}
\pubyear{2014}
\firstpage{2039}
\lastpage{2075}
\doi{10.3150/13-BEJ550} %kopijuoti is PTS

\makeatletter
\newcommand{\rrVert}{\Vert}
\newcommand{\rrvert}{\vert}
\newcommand{\llVert}{\Vert}
\newcommand{\llvert}{\vert}
\newcommand{\trup}[2]{{#1}/{#2}}
\renewcommand{\emptyset}{\varnothing}
\newcommand{\eqref}[1]{(\ref{#1})}
\newcommand{\Pb}{\mathbb{P}}
\newcommand{\R}{\mathbb{R}}
\newcommand{\N}{\mathbb{N}}
\newcommand{\one}{\mathbh{1}}

\newtheorem{theorem}{Theorem}
\newtheorem{lemma}{Lemma}
\newtheorem{prop}{Proposition}
\newproclaim{Assumption}{Assumption}
\makeatother

\begin{document}
\begin{frontmatter}

\title{About the posterior distribution in hidden Markov models with
unknown number of states}
\runtitle{Posterior distribution for HMM}

\begin{aug}
%%%% inicialai - be tarpu
\author[1]{\inits{E.}\fnms{Elisabeth} \snm{Gassiat}\corref{}\thanksref{1}\ead[label=e1]{elisabeth.gassiat@math.u-psud.fr}} \and
\author[2]{\inits{J.}\fnms{Judith} \snm{Rousseau}\thanksref{2}\ead[label=e2]{rousseau@ceremade.dauphine.fr}}
%%\runauthor{} %% auto
\address[1]{Laboratoire de Math\'ematiques d'Orsay UMR 8628,
Universit\'e Paris-Sud,
B\^atiment 425, 91405 Orsay-C\'edex, France. \printead{e1}}
\address[2]{CREST-ENSAE,
3 avenue Pierre Larousse,
92245 Malakoff Cedex, France.\\ \printead{e2}}
\end{aug}

% HISTORY:
\received{\smonth{3} \syear{2013}}
\revised{\smonth{7} \syear{2013}}

% ABSTRACT
%
\begin{abstract}
We consider finite state space stationary hidden Markov models (HMMs)
in the situation where the number of hidden states is unknown. We
provide a frequentist asymptotic evaluation of Bayesian analysis
methods. Our main result gives posterior concentration rates for the
marginal densities, that is for the density of a fixed number of
consecutive observations. Using conditions on the prior, we are then
able to define a consistent Bayesian estimator of the number of hidden states.
It is known that the likelihood ratio test statistic for overfitted
HMMs has a nonstandard behaviour and is unbounded. Our conditions on
the prior may be seen as a way to penalize parameters to avoid this phenomenon.
Inference of parameters is a much more difficult task than inference of
marginal densities,
we still provide a precise description of the situation when the
observations are i.i.d. and we allow for $2$ possible hidden states.
\end{abstract}

% KEYWORDS
% visi is mazosios raides ir pagal abecele
%

\end{frontmatter}

%s1 #&#
\section{Introduction}
\label{sec:intro}

Finite state space hidden Markov models (which will be shortened to
HMMs throughout the paper) are stochastic processes $(X_j,Y_j)_{j\geq
1}$ where $(X_j)_{j\geq1} $ is a Markov chain living in a finite state
space $\mathcal X$ and conditionally on $(X_j)_{j\geq1}$ the $Y_j$'s
are independent with a distribution depending only on $X_j$ and living
in $\mathcal Y$.
%The observations are $Y_{1:n} = (Y_1,\cdots,Y_n)$ and the associated
%states $X_{1:n} = (X_1,\cdots,X_n)$ are unobserved.
HMMs are useful tools to model time series where the observed
phenomenon is driven by a latent Markov chain.
They have been used successfully in a variety of applications,
% such as economics (e.g. \citet{albert:chib:1993a}), genomics (e.g.
%speech recognition (e.g. \citet{rabiner:89}) to name but a few. T
the books MacDonald and Zucchini \cite{macdonald:zucchini:1997},
Zucchini and MacDonald \cite{macdonald:zucchini:2009} and Capp{\'e}
\textit{et al.} \cite{cappe:moulines:ryden:2004} provide
several examples of applications of HMMs and give a recent (for the
latter) state of the art in the statistical analysis of HMMs.
Finite state space HMMs may also be seen as a dynamic extension of
finite mixture models and may be used to do unsupervised clustering.
%The number of hidden states induces a classification of the regimes in
%which the time series evolves.
The hidden states often have a practical interpretation in the
modelling of the underlying phenomenon. It is thus of importance to be
able to infer both the number of hidden states (which we call the order
of the HMM) from the data, and the associated parameters.

The aim of this paper is to provide a frequentist asymptotic analysis
of Bayesian methods used for statistical inference in finite state
space HMMs when the order is unknown. Let us first review what is known
on the subject and important questions that still stay unsolved.

In the frequentist literature, penalized likelihood methods have been
proposed to estimate the order of a HMM, using for instance Bayesian
information criteria (BIC for short). These methods were applied for
instance in Leroux and Putterman \cite{leroux:putterman:1992}, Ryd{\'e}n
\textit{et al.} \cite{rydenetal:1998}, but
without theoretical consistency results. Later, it has been observed
that the likelihood ratio statistics is unbounded, in the very simple
situation where one wants to test between $1$ or $2$ hidden states,
see Gassiat and K\'eribin \cite{CKEG}. The question whether BIC penalized
likelihood methods
lead to consistent order estimation stayed open. Using tools borrowed
from information theory, it has been possible to calibrate heavier
penalties in maximum likelihood methods to obtain consistent estimators
of the order, see Gassiat and Boucheron \cite{GB03}, Chambaz \textit{et al.}
\cite{CGG09}.
The use of penalized marginal pseudo likelihood was also proved to lead
to weakly consistent estimators by Gassiat \cite{Gas02}.

On the Bayesian side, various methods were proposed to deal with an
unknown number of hidden states, but no frequentist theoretical result
exists for these methods.
Notice though that, if the number of states is known, de~Gunst and
Shcherbakova \cite{degunst:shcherbakova:08} obtain a Bernstein--von Mises theorem for
the posterior distribution, under additional (but usual) regularity conditions.
When the order is unknown, reversible jump methods have been built,
leading to satisfactory results on simulation and real data, see Boys and
Henderson \cite{boys:henderson:2004}, Green and Richardson
\cite{green:richardson:2002}, Robert \textit{et al.}
\cite{robert:ryden:titterington:1999b}, Spezia \cite{Spezia09}. The ideas of
variational Bayesian methods were developed in
McGrory and Titterington \cite{titterington:mcgrory:2009}.
Recently, one of the authors proposed a frequentist asymptotic analysis
of the posterior distribution for overfitted mixtures when the
observations are i.i.d., see Rousseau and Mengersen \cite{kerrieju:11}.
In this paper, it is
proved that one may choose the prior in such a way that extra
components are emptied, or in such a way that extra components merge
with true ones. More precisely, if a Dirichlet prior $\mathcal {D}(\alpha
_1,\ldots,\alpha_k)$ is considered on the $k$ weights of the mixture
components, small values of the $\alpha_j$'s imply that the posterior
distribution will tend to empty the extra components of the mixture
when the true distribution has a smaller number, say $k_0 < k$ of true
components.
%So, this paper gives a guideline for the choice of the prior so as to
%be able to interpret the posterior. \\
One aim of our paper is to understand if such an analysis may be
extended to HMMs.

As is well known in the statistical analysis of overfitted finite
mixtures, the difficulty of the problem comes from the
non-identifiability of the parameters.
But what is specific to HMMs is that the non-identifiability of the
parameters leads to the fact that neighbourhoods of the ``true''
parameter values contain transition matrices arbitrarily close to non-ergodic transition matrices.
To understand this on a simple example, just consider the case of HMMs
with two hidden states, say $p$ is the probability of going from state
1 to state 2 and $q$ the probability of going from state 2 to state 1.
If the observations are in fact independently distributed, their
distribution may be seen as a HMM with two hidden states where $q=1-p$.
Neighbourhoods of the ``true'' values $(p,1-p)$ contain parameters such
that $p$ is small or $1-p$ is small, leading to hidden Markov chains
having mixing coefficients very close to $1$. Imposing a prior
condition such as $\delta\leq p \leq1-\delta$ for some $\delta>0$ is
not satisfactory.

Our first main result Theorem~\ref{theo:gene:finite} gives
concentration rates for the posterior distribution of the marginal
densities of a fixed number of consecutive observations. First, under
mild assumptions on the densities and the prior, we obtain the
asymptotic posterior concentration rate $\sqrt{n}$, $n$ the number of
observations, up to a $\log n$ factor, when the loss function is the
$L_1$ norm between densities multiplied by some function of the
ergodicity coefficient of the hidden Markov chain. Then, with more
stringent assumptions on the prior, we give posterior concentration
rates for the marginal densities in $L_1$ norm only (without the
ergodicity coefficient). For instance,
consider a finite state space HMM, with $k$ states and with independent
Dirichlet prior distributions $\mathcal{D}(\alpha_1,\ldots,\alpha_k)$ on
each row of the transition matrix of the latent Markov chain. Then our
theorem says that if the sum of the parameters $\alpha_j$'s is large
enough, the posterior distribution of the marginal densities in $L_1$
norm concentrates at a polynomial rate in $n$. These results are
obtained as applications of a general theorem we prove about
concentration rates for the posterior distribution of the marginal
densities when the state space of the HMM is not constrained to be a
finite set, see Theorem~\ref{theo:gene}.

A byproduct of the non-identifiability for overfitted mixtures or HMMs
is the fact that, going back from marginal densities to the parameters
is not easy.
The local geometry of finite mixtures has been understood by Gassiat and
van Handel \cite{HanGas}, and following their approach in the HMM context
we can go
back from the $L_1$ norm between densities to the parameters.
We are then able to propose a Bayesian consistent estimator of the
number of hidden states, see Theorem~\ref{theo:order}, under the same
conditions on the prior as in Theorem~\ref{theo:gene:finite}. To our
knowledge, this is the first consistency result on Bayesian order
estimation in the case of HMMs.

Finally, obtaining posterior concentration rates for the parameters
themselves seems to be very difficult, and we propose a more complete
analysis in the simple situation of HMMs with $2$ hidden states and
independent observations. In such a case, we prove that, if all the
parameters (not only the sum of them) of the prior Dirichlet
distribution are large enough, then extra components merge with true
ones, see Theorem~\ref{theo:uncontredeux}. We believe this to be more
general but have not been able to prove it.

The organization of the paper is the following. In Section~\ref{sec:finite},
we first set the model and notations. In subsequent
subsections, we give Theorems~\ref{theo:gene:finite},
\ref{theo:order} and~\ref{theo:uncontredeux}. In Section~\ref{sec:gene},
we give the posterior concentration theorem for general
HMMs, Theorem~\ref{theo:gene}, on which Theorem~\ref{theo:gene:finite}
is based. All proofs are given in Section~\ref{sec:proofs}.

%%%%%%%%%%%%%%%%%%%%%%%%%%%
%s2 #&#
\section{Finite state space hidden Markov models}
\label{sec:finite}

%s2.1 #&#
\subsection{Model and notations}
\label{subsec:notation}

Recall that finite state space HMMs model pairs $(X_i, Y_i)_{i\geq1}$
where $(X_i)_{i\geq1}$ is the unobserved Markov chain living on a
finite state space $\mathcal X = \{1, \ldots, k\}$ and the observations
$(Y_i)_{i\geq1}$ are conditionally independent given the $(X_i)_{i\geq
1}$. The observations take value in $\mathcal Y$,
%. The space $\mathcal Y$
which is assumed to be a Polish space endowed with its $\sigma$-field.
Throughout the paper, we denote $x_{1:n} = (x_1, \ldots, x_n)$.

The hidden Markov chain $(X_i)_{i\geq1}$ has a Markov transition
matrix $Q = (q_{ij})_{1\leq i, j \leq k}$. The conditional distribution
of $Y_i $ given $X_i$ has a density with respect to some given measure
$\nu$ on $\mathcal Y$. We denote by $g_{\gamma_j}(y)$, $j= 1,\ldots,k$,
the conditional density of $Y_i $ given $X_i =j$. Here, $\gamma_j \in
\Gamma\subset\R^d$ for $j =1,  \ldots, k$, the $\gamma_{j}$'s are
called the emission parameters. In the following, we parametrize the
transition matrices on $\{1,\ldots,k\}$ as $(q_{ij})_{1\leq i \leq k,
1\leq j \leq k-1}$ (implying that $q_{i k}=1-\sum_{j=1}^{k-1}q_{ij}$
for all $i\leq k$)
and we denote by $\Delta_k$ the set of probability mass functions
$ \Delta_{k}=\{(u_{1},\ldots, u_{k-1}) \dvtx  u_{1}\geq0,\ldots,
u_{k-1}\geq0,\sum_{i=1}^{k-1}u_{i} \leq1 \}$.
We shall also use the set of positive probability mass functions
$
\Delta_{k}^{0}=\{(u_{1},\ldots, u_{k-1}) \dvtx  u_{1}> 0,\ldots, u_{k-1}>
0,\sum_{i=1}^{k-1}u_{i} < 1 \}
$.
Thus, we may denote the overall parameter by $\theta= (q_{ij}, 1\leq i
\leq k, 1\leq j \leq k-1 ; \gamma_{1},\ldots,\gamma_{k} ) \in\Theta
_{k}$ where
$\Theta_{k}= {\Delta}_{k}^{k} \times\Gamma^{k}$. To alleviate
notations, we will write $\theta= (Q ; \gamma_{1},\ldots,\gamma_{k} )$,
where $Q=(q_{ij})_{ 1\leq i, j \leq k}$, $q_{ik}=1-\sum_{j=1}^{k-1}q_{ij}$ for all $i\leq k$.

Throughout the paper, $ \nabla_{\theta}h$ denotes the gradient vector
of the function $h$ when considered as a function of $\theta$, and
$D^{i}_{\theta}h$ its $i$th derivative operator with respect to $\theta
$, for $i \geq1$.
We denote by $B_d(\gamma, \epsilon)$ the $d$ dimensional ball centered
at $\gamma$ with radius $\epsilon$, when $\gamma\in\R^d $. The
notation $a_{n}\gtrsim b_{n}$ means that $a_{n}$ is larger than $b_{n}$
up to a positive constant that is fixed throughout.

%We denote by $D^{i}_{\gamma} g_{\gamma_1}$ the $i$-th derivative of $g_
%$i \geq0$.

Any Markov chain on a finite state space with transition matrix $Q$
admits a stationary distribution which we denote by $\mu_Q$, if it
admits more than one we choose one of them. Then for any finite state
space Markov chain with transition matrix $Q$
it is possible to define real numbers $\rho_{Q}\geq1$
such that, for any integer $m$, any $j \leq k$
%
%e1 #&#
\begin{equation}
 \sum_{j=1}^k\bigl\llvert
\bigl(Q^m\bigr)_{ij} - \mu_Q(j)\bigr\rrvert
\leq \rho_{Q}^{-m},\qquad \rho_Q = \Biggl( 1 -
\sum_{j=1}^k \min_{1 \leq i \leq k}
q_{ij} \Biggr)^{-1},
\end{equation}
where $Q^m$ is the $m$-step transition matrix of the Markov chain.
If $\rho_{Q}>1$, the Markov chain $(X_{n})_{n\geq1}$ is uniformly
geometrically ergodic and $\mu_{Q}$ is its unique stationary
distribution. In the following, we shall also denote $\mu_\theta$ and
$\rho_\theta$ in the place of $\mu_Q$ and $\rho_Q$ when $\theta= (Q ;
\gamma_1, \ldots, \gamma_k)$.

We write $\Pb_{\theta}$ for the probability distribution of the
stationary HMM $(X_{j},Y_{j})_{j\geq1}$ with parameter $\theta$. That
is, for any integer $n$, any set $A$ in the Borel $\sigma$-field of
${\mathcal X}^{n}\times{\mathcal Y}^{n}$:
%
%e2 #&#
\begin{eqnarray}
\label{defhmmstat} %
&& \Pb_{\theta} \bigl((X_{1},
\ldots,X_{n},Y_{1},\ldots,Y_{n})\in A \bigr)
\nonumber
\\[-8pt]
\\[-8pt]
&&\quad =\sum_{x_1,\ldots, x_n=1}^k \int
_{\mathcal Y^n} \one _{A}(x_{1:n},y_{1:n})
\mu_{Q}(x_{1}) \prod_{i=1}^{n-1}q_{x_ix_{i+1}}
\prod_{i=1}^{n}g_{\gamma_{x_i} }
(y_{i} )\nu(\mathrm{d}y_{1})\cdots\nu(\mathrm{d}y_{n}).
\nonumber
\end{eqnarray}
Thus for any integer $n$, under $\Pb_{\theta}$, $Y_{1:n}=(Y_{1},\ldots
,Y_{n})$ has a probability density with respect to $\nu(\mathrm{d}y_{1})\cdots\nu
(\mathrm{d}y_{n})$ equal to
%
%e3 #&#
\begin{equation}
\label{margeq} f_{n,\theta}(y_1,\ldots, y_n) = \sum
_{x_1,\ldots, x_n=1}^k \mu_{Q}(x_{1})
\prod_{i=1}^{n-1}q_{x_ix_{i+1}}\prod
_{i=1}^{n}g_{\gamma_{x_i} }
(y_{i} ).
\end{equation}
We note $E_{\theta}$ for the expectation under $\Pb_{\theta}$.

We denote $\Pi_k$ the prior distribution on $\Theta_k$. As is often the
case in Bayesian analysis of HMMs, instead of computing the stationary
distribution $\mu_Q$ of the hidden Markov chain with transition matrix
$Q$, we consider a probability distribution $\pi_{\mathcal X}$ on the
unobserved initial state $X_0$. Denote $\ell_n(\theta,x_0)$ the
log-likelihood starting from $x_0$, for all $x_0 \in\{1,\ldots, k\}$,
we have
\[
\ell_{n} (\theta,x_0 )=\log \Biggl[\sum
_{x_1,\ldots, x_n=1}^k \prod_{i=0}^{n-1}q_{x_ix_{i+1}}
\prod_{i=1}^{n}g_{\gamma_{x_i} }
(y_{i} ) \Biggr].
\]
The log-likelihood starting from a probability distribution $ \pi
_{\mathcal X}$ on $\mathcal X$ is then given by
%denoted
%$\ell_{n}\left(\theta,\pi_0\right)$ i.e.
$
%$\ell_{n}\left(\theta,\pi_0\right) =
\log [ \sum_{x_0=1}^k \mathrm{e}^{\ell_{n} (\theta,x_0 )} \pi
_{\mathcal X}(x_0 )  ]$.
%$
This may\vadjust{\goodbreak} also be interpreted as taking a prior $ \Pi= \Pi_k\otimes\pi
_{\mathcal X}$ over $\Theta_k \times\{ 1, \ldots, k\}$.
The posterior distribution can then be written as
%
%e4 #&#
\begin{equation}
\label{posterior} \Pb^{ \Pi} (A \vert Y_{1:n} ) =
\frac{ \sum_{x_0=1}^k\int_{A}
\mathrm{e}^{\ell_n(\theta,x_0)} \Pi_k (\mathrm{d}\theta ) \pi_{\mathcal
X}(x_0)}{ \sum_{x_0=1}^k\int_{\Theta} \mathrm{e}^{\ell_n(\theta,x_0)} \Pi_k
(\mathrm{d}\theta ) \pi_{\mathcal X}(x_0) }
\end{equation}
for any Borel set $A\subset\Theta_k$.

%We shall also use the notation $\Pb_{\theta,x_0}$ for the probability
%distribution of the HMM starting from $x_0$, that is, for any integer
%$n$, any measurable set $A$ in the Borel $\sigma$-field of ${\mathcal
%X}^{n}\times{\mathcal Y}^{n}$:
%so that for any $\theta\in\Theta_k$,
%$$
%$$
%We denote by $E_{\theta}$ the expectation under $\Pb_{\theta}$ and by
%$E_{\theta,x}$ the expectation under $\Pb_{\theta,x}$.\\

Let $\mathcal{M}_{k}$ be the set of all possible probability
distributions
%of $(Y_{n})_{n\geq1}$ under
$\Pb_{\theta}$ for all $\theta\in\Theta_k$. We say that the HMM $\Pb
_{\theta}$ has order $k_{0}$ if the probability distribution of
$(Y_{n})_{n\geq1}$ under $\Pb_{\theta}$ is in
$\mathcal{M}_{k_{0}}$ and not in $\mathcal{M}_{k}$ for all $k<k_{0}$.
Notice that a HMM of order $k_{0}$ may be represented as a HMM of order
$k$ for any $k>k_{0}$. Indeed, let $Q^{0}$ be a $k_{0}\times k_{0}$
transition matrix, and $(\gamma_{1}^{0},\ldots,\gamma_{k_0}^{0})\in
\Gamma^{k_0}$ be parameters that define a HMM of order $k_{0}$. Then,
$\theta= (Q ; \gamma_{1}^{0},\ldots,\gamma_{k_0}^{0},\ldots,\gamma
_{k_0}^{0} )\in\Theta_{k}$ with $Q=(q_{ij}, 1\leq i,j \leq k)$ such that:
%
%e5 #&#
\begin{equation}
\label{lem:m1} %
\begin{array} {@{}rcl@{\qquad}l@{}} q_{ij } &
=& q_{ij}^{0},& i, j < k_{0},
\\[9pt]
q_{ij} & =& q_{k_{0}j}^{0},& i\geq
k_{0}, j<k_{0},
\\[9pt]
\displaystyle\sum_{l=k_0}^{k}q_{il}
&= &q_{ik_0}^{0},& i\leq k_0, \quad\mbox{and}
\quad \displaystyle\sum_{l=k_0}^{k}q_{il}=q_{k_0k_0}^{0},
\qquad i\geq k_0 \end{array} %
\end{equation}
gives $\Pb_{\theta}=\Pb_{\theta_0}$.
Indeed, let $(X_{n})_{n\geq1}$ be a Markov chain on $\{1,\ldots,k\}$
with transition matrix $Q$. Let $Z$ be the function from $\{1,\ldots,k\}
$ to $\{1,\ldots,k_0\}$ defined by $Z(x)=x$ if $x\leq k_0$ and
$Z(x)=k_0$ if $x\geq k_0$. Then $(Z(X_{n}))_{n\geq1}$ is a Markov
chain on $\{1,\ldots,k_0\}$ with transition matrix $Q^0$.

%%%%%%%%%%%%%%%%%%%%%%%%

%s2.2 #&#
\subsection{Posterior convergence rates for the finite marginal densities}
\label{subsec:finite:gene}

Let $\theta_0 = (Q^{0} ; \gamma_{1}^{0},\ldots,\gamma_{k_0}^{0})\in
\Theta_{k_0}$, $Q^{0}=(q_{ij}^{0})_{ 1\leq i \leq k_0, 1\leq j \leq
k_0}$, be the parameter of a HMM of order $k_0 \leq k$. We now assume
that $\Pb_{\theta_{0}}$ is the distribution of the observations.
In this section, we fix an integer $l$ and study the posterior
distribution of the density of $l$ consecutive observations, that is
$f_{l,\theta}$, given by \eqref{margeq}
with $n=l$.
We
study the posterior concentration rate around $f_{l,  \theta_0}$ in
terms of the $L_1$ loss function, when $\Pb_{\theta_0}$ is possibly
of order $k_{0}<k$.
%non identifiable, i.e. the Markov chain lives on a state space with
%fewer number of states.
In this case, Theorem~2.1 of de~Gunst and Shcherbakova \cite{degunst:shcherbakova:08} does not
apply and there is no result in the literature
% describing the posterior concentration rate.
about the frequentist asymptotic properties of the posterior distribution.
The interesting and difficult feature of this case is that even though
$\theta_0$ is parameterized as an ergodic Markov chain $Q^0$ with $k$
states and some identical emission parameters as described in
\eqref{lem:m1}, $f_{l,  \theta_0}$ can be approached by marginals $f_{l,
\theta}$ for which $\rho_\theta$ is arbitrarily close to 1, which
deteriorates the posterior concentration rate, see Theorem~\ref{theo:gene:finite}.
%
% Interestingly,
%in the case of HMM, the identifiability issue described above has not
%only an impact on the study of the asymptotic behaviour of the
%posterior distribution in $\theta$ but it also has an impact on the
%study of the asymptotic posterior distribution of the marginal
%densities. The $l$-dimensional marginal density, for $l \geq1$, is
%given by \eqref{margeq}
%f_{l,\theta}(y)=\sum_{1\leq i_{1},\ldots,i_{l} \leq k}\mu_{\theta}
%with $\theta= (q_{ij}, 1\leq i \leq k, 1\leq j \leq k-1 ;
%in ${\mathcal Y}^{l}$.
%with $n=l$.
%this is in sharp contrast with the case of independent mixtures, i.e.
%when the $X_i$ are independent and identically distributed, for which
%the posterior distribution on $f_{1,\theta}$ concentrates around
%$f_{1,\theta_0}$ at the rate $1/\sqrt{n}$ (possibly up to a $\log n$
%term) in the $L_1$ sense; see for instance \citet{kerrieju:11}.

Let $\pi(u_{1},\ldots,u_{k-1}) $ be a prior density with respect to
the Lebesgue measure on $\Delta_k$, and let $\omega(\gamma) $ be a
prior density on $\Gamma$ (with respect to the Lebesgue measure on $\R
^{d}$). We consider prior distributions such that the rows of the
transitions matrix $Q$ are independently distributed from $\pi$ and
independent of the component parameters $\gamma_i$, $i=1,\ldots,k$,
which are independently distributed from $\omega$. Hence, the prior
density of $\Pi_k$ (with respect to the Lebesgue measure) is equal to
$\pi_{k} = \pi^{\otimes k} \otimes\omega^{\otimes k}$. We still
denote by $\pi_{\mathcal X}$ a probability on $\{1,\ldots, k\}$, we
assume that $\pi_{\mathcal X}(x) >0$ for all $x\in\{1,\ldots, k\}$ and
set $\Pi= \Pi_k \otimes\pi_{\mathcal X}$.
% In this section we use a
%Dirichlet type prior, see assumption \textbf{A3} below, or an
%exponential type prior, see assumption (FE1) below, on the transition
%parameters $(q_{ij}, j\leq k)$.
We shall use the following assumptions.
\begin{enumerate}[\textbf{A3}]
\item[\textbf{A0}]
$q_{ij}^{0} >0$, $ 1\leq i \leq k_0$, $1\leq j \leq k_0 $.\vspace*{2pt}
\item[\textbf{A1}]
The function $\gamma\mapsto g_{\gamma}(y)$ is twice continuously
differentiable in $\Gamma$, and for any $\gamma\in\Gamma$, there exists
$\epsilon>0$ such that
\begin{eqnarray*}
\int\sup_{\gamma'\in B_d(\gamma, \epsilon)}\bigl \| \nabla_{\gamma} \log
g_{\gamma'} (y )\bigr \|^{2}g_{\gamma} (y )\nu(\mathrm{d}y)&<&
+\infty,
\\
\int\sup_{\gamma'\in B_d(\gamma, \epsilon)}\bigl \| D^{2}_{\gamma} \log
g_{\gamma'} (y )\bigr \|^{2}g_{\gamma} (y )\nu(\mathrm{d}y) &<&
+\infty,
\end{eqnarray*}
$\|\sup_{\gamma'\in B_d(\gamma, \epsilon)} \nabla_{\gamma} g_{\gamma'}
(y )\| \in L_{1}(\nu)$ and $\| \sup_{\gamma'\in B_d(\gamma,
\epsilon)} D^{2}_{\gamma} g_{\gamma'}  (y )\| \in L_{1}(\nu)$.
\item[\textbf{A2}]
There exist $a>0$ and $b>0$ such that
\[
\sup_{\| \gamma\| \leq n^{b}}\int \bigl \| \nabla_{\gamma} g_{\gamma} (y
)\bigr \| \,\mathrm{d}\nu(y) \leq n^{a}.
\]
\item[\textbf{A3}]
$\pi$ is continuous and positive on $\Delta_{k}^{0}$, and there exists
$C, \alpha_{1}>0, \ldots, \alpha_{k}>0$ such that (Dirichlet type priors):
\begin{eqnarray*}
\forall(u_{1},\ldots,u_{k-1})\in\Delta_{k}^{0},
\qquad u_k &=& 1- \sum_{i=1}^{k-1}u_i,
\\
 0&<&\pi (u_{1},\ldots,u_{k-1} ) \leq C
u_{1}^{\alpha_{1}-1}\cdots u_k^{\alpha_{k}-1}
\end{eqnarray*}
and $\omega$ is continuous and positive on $\Gamma$ and satisfies
%
%e6 #&#
\begin{equation}
\label{cond:w} \int_{\|x\|\geq n^b } \omega(x)\,\mathrm{d}x =
\mathrm{o}\bigl(n^{-k(k-1+d)/2}\bigr),
\end{equation}
with $b$ defined in assumption \textbf{A2}.
\end{enumerate}
We will alternatively replace \textbf{A3} by
\begin{enumerate}[\textbf{A3bis}]
\item[\textbf{A3bis}]
$\pi$ is continuous and positive on $\Delta_{k}^{0}$, and there exists
$C$ such that (exponential type priors):
\begin{eqnarray*}
\forall(u_{1},\ldots,u_{k-1})\in\Delta_{k}^{0},
\qquad u_k &=& 1- \sum_{i=1}^{k-1}u_i,
\\
0&<& \pi (u_{1},\ldots,u_{k-1} ) \leq C
\exp(-C/u_{1} )\cdots \exp(-C/u_{k} )
\end{eqnarray*}
and $\omega$ is continuous and positive on $\Gamma$ and satisfies \eqref{cond:w}.
\end{enumerate}

%th1 #&#
\begin{theorem}
\label{theo:gene:finite}
Assume \textup{\textbf{A0}--\textbf{A3}}.
%, \textbf{A3}, \textbf{A1} and \textbf{A2}.
Then, there exists $K$ large enough such that
%
%e7 #&#
\begin{equation}
\label{conc:rho} \Pb^{\Pi} \biggl[   \theta\dvtx \|f_{l,\theta}
- f_{l,\theta_0}\|_1 (\rho_\theta- 1) \geq K \sqrt{
\frac{\log n}{n}} \Big\rrvert Y_{1:n} \biggr] = \mathrm{o}_{\Pb_{\theta_0}}
(1 ),
\end{equation}
where $\rho_{\theta}= (1-\sum_{j=1}^{k}\inf_{1\leq i \leq
k}q_{ij} )^{-1}$. If moreover $\bar{\alpha} := \sum_{1\leq i \leq
k}\alpha_{i} > k(k-1 + d)$, then
%
%e8 #&#
\begin{equation}
\label{fl:conc:A3} \Pb^{\Pi} \bigl[   \theta\dvtx \|f_{l,\theta}
- f_{l,\theta_0}\|_1 \geq2K n ^{-\trup{(\bar{\alpha}-k( k-1 + d) )}
{( 2\bar{\alpha} )} } (\log n ) \big\rrvert
Y_{1:n} \bigr] = \mathrm{o}_{\Pb_{\theta_0}} (1 ).
\end{equation}
If we replace \textup{\textbf{A3}} by \textup{\textbf{A3bis}},
then there exists $K$ large enough such that
%
%e9 #&#
\begin{equation}
\label{fl:conc:A3bis} \Pb^{\Pi} \bigl[  \theta\dvtx \|f_{l,\theta} -
f_{l,\theta_0}\|_1 \geq 2K n ^{-1/2}(\log n
)^{3/2} \big\rrvert Y_{1:n} \bigr] = \mathrm{o}_{\Pb_{\theta
_0}} (1
).
\end{equation}
\end{theorem}

Theorem~\ref{theo:gene:finite} is proved in Section~\ref{subsec:gene:finite}
as a consequence of Theorem~\ref{theo:gene} stated
in Section~\ref{sec:gene}, which gives posterior concentration rates
for general HMMs.

Assumption \textbf{A0} is the usual ergodic condition on the finite
state space Markov chain.
Assumptions \textbf{A1} and \textbf{A2} are mild usual regularity
conditions on the emission densities $g_\gamma$ and hold for instance
for multidimensional Gaussian distributions, Poisson distributions, or
any regular exponential families. Assumption \textbf{A3} on the prior
distribution of the transition matrix $Q$ is satisfied for instance if
each row of $Q$ follows a Dirichlet distribution or a mixture of
Dirichlet distributions, as used in Nur \textit{et al.} \cite{nur:etal:09},
and assumption
\eqref{cond:w} is verified for densities $\omega$ that have at most
polynomial tails.
%The lightness of the tails depend on the constraints on $b$. In many
%instances of parametric families $b$ can be chosen arbitrarily so that
%any polynomial tail on $\omega$ verifies \eqref{cond:w}. Condition
%is close to a truncated prior on an open subset of $\Delta_{k}^0$.

The constraint on $\bar\alpha= \sum_{i}\alpha_i$ or condition
\textbf{A3bis} are used to ensure that \eqref{fl:conc:A3} and
\eqref{fl:conc:A3bis}
hold respectively. The posterior concentration result \eqref{conc:rho}
implies that the posterior distribution might put non-negligible mass
on values of $\theta$ for which $\rho_\theta-1$ is small and $\|
f_{l,\theta} - f_{l,\theta_0}\|_1$ is not. These are parameter values
associated to nearly non-ergodic latent Markov chains. Since $\rho
_\theta-1$ is small is equivalent to $\sum_{j} \min_i q_{ij}$ is small,
the condition $\bar\alpha> k(k-1+d)$ prevents such pathological
behaviour by ensuring that the prior mass of such sets is small enough.
This condition is therefore of a different nature than Rousseau and
Mengersen's \cite{kerrieju:11} condition on the prior, which characterizes the
asymptotic behaviour of the posterior distribution on the parameter
$\theta$. In other words, their condition allows in (static) mixture
models to go from a posterior concentration result on $f_{l,\theta}$ to
a posterior concentration result on $\theta$ whereas, here, the
constraint on $\bar{\alpha}$ is used to obtain a posterior concentration
result on $f_{l,\theta}$.
Going back from $\|f_{l,\theta} - f_{l,\theta_0}\|_1$ to the parameters
requires a deeper understanding of the geometry of finite HMMs, similar
to the one developed in Gassiat and van Handel \cite{HanGas}. This will be needed to estimate
the order of the HMM in Section~\ref{subset:number}, and fully explored
when $k_{0}=1$ and $k=2$ in Section~\ref{subsec:1et2}.

For general priors, we do not know whether the $\sqrt{\log n}$ factor
appearing in \eqref{conc:rho} could be replaced or not by any sequence
tending to infinity.
In the case where the $\alpha_i$'s are large enough (Dirichlet type
priors), and when $k_0=1$ and $k=2$, we obtain a concentration rate
without the $\sqrt{\log n}$ factor, see Lemma~\ref{lemrootn} in
Section~\ref{theo:uncontredeux}. To do so, we prove Lemma~\ref{k2k01:Dn} in
Section~\ref{theo:uncontredeux} for which we need to compute explicitly
the stationary distribution and the predictive probabilities to obtain
a precise control of the likelihood, for $\theta$'s such that $\Pb
_{\theta}$ is near $\Pb_{\theta_0}$, and to control local entropies of
slices for $\theta$'s such that $\Pb_{\theta}$ is near $\Pb_{\theta_0}$
and where $\rho_{\theta}-1$ might be small. It is not clear to us that
extending such computations to the general case is possible in a
similar fashion. The $\log n$ terms appearing in \eqref{fl:conc:A3} and
\eqref{fl:conc:A3bis} are consequences of the $\sqrt{\log n}$ term
appearing in \eqref{conc:rho}.

%%%%%%%%%%%%%%%

%s2.3 #&#
\subsection{Consistent Bayesian estimation of the number of states}
\label{subset:number}

To define a Bayesian estimator of the number of hidden states $k_0$, we
need to decide how many states have enough probability mass, and are
such that their emission parameters are different enough.
We will be able to do it under the assumptions of Theorem~\ref{theo:gene:finite}.
Set
$w_{n}=n^{-(\bar\alpha- k(k+d-1))/(2\bar\alpha) }\log n
$ if \textbf{A3} holds and $\bar\alpha> k(k+d-1) $, and set
$w_{n}=n^{-1/2 }(\log n)^{3/2}
$ if instead \textbf{A3bis} holds.
Let $(u_{n})_{n\geq1}$
and $(v_{n})_{n\geq1}$
be sequences of positive real numbers tending to $0$ as $n$ tends to
infinity such that $w_{n}=\mathrm{o}(u_{n}v_{n})$. As in Rousseau and
Mengersen \cite{kerrieju:11}, in
the case of a misspecified model with $k_0 <k$, $f_{l, \theta_0}$ can
be represented by merging components or by emptying extra components.
For any $\theta\in\Theta_{k}$, we thus define $J (\theta )$ as
\[
J (\theta )= \bigl\{j \dvtx \Pb_{\theta} (X_{1}=j )\geq
u_{n} \bigr\}, %u_{n}^{1/3}
\]
that is, $J(\theta)$ corresponds to the set of \textit{non-empty}
components. To cluster the components that have similar emission
parameters, we define for all $j \in J(\theta)$
%$$A_{j}\left(\theta\right)=\left\{i\in J\left(\theta\right) : \Pb_{
%u_{n}v_{n} \right\}.
%$$
\[
A_{j} (\theta )= \bigl\{i\in J (\theta ) \dvtx \|\gamma
_{j}-\gamma_{i}\|^2 \leq v_n \bigr
\} %\frac{w_{n}}{u_{n}v_{n}}
\]
and the clusters are defined by: for all $j_{1}, j_2 \in J(\theta)$,
$j_1$ and $j_2$ belong to the same cluster (noted $j_1 \sim j_2$) if
and only if there exist $r>1$ and $i_1, \ldots, i_r\in J(\theta)$ with
$i_1 = j_1$ and $i_r = j_2$ such that for all $1 \leq l \leq r-1$,
$A_{i_l} (\theta )\cap A_{i_{l+1}} (\theta )\neq
\emptyset$.
We then
define the effective order of the HMM at $\theta$ as the number $L
(\theta )$ of different clusters, that is, as the number of
equivalent classes with respect to the equivalence relation $\sim$
defined above. By a good choice of $u_n$ and $v_n$, we construct a
consistent estimator of $k_0$ by considering either the posterior mode
of $L(\theta)$ or its posterior median. This is presented in Theorem~\ref{theo:order}.

To prove that this gives a consistent estimator, we need an inequality
that relates the $L_1$ distance between the $l$-marginals, $\|
f_{l,\theta} - f_{l,  \theta_0}\|_1$, to a distance between the
parameter $\theta$ and parameters $\tilde{\theta}_{0}$ in $\Theta_{k}$
such that $f_{l,\tilde{\theta}_0}=f_{l,\theta_0}$.
% $\theta$ and some representation of $\theta_0$ (recall that $f_{l,
Such an inequality will be proved in Section~\ref{subsec:order}, under
the following structural assumption.

%To define the clustering structure,
Let $T=\{\mathbf{t}=(t_{1},\ldots,t_{k_{0}})\in\{1,\ldots,k\}^{k_{0}}\dvtx
t_{i}<t_{i+1}, i=0,\ldots,k_{0}-1\}$. If $b$ is a vector, $b^T$ denotes
its transpose.
\begin{enumerate}[\textbf{A4}]
\item[\textbf{A4}]
For any $\mathbf{t}=(t_{1},\ldots,t_{k_{0}})\in T$,
any $(\pi_i)_{i = 1}^{k - t_{k_0}} \in(\R^+)^{ k - t_{k_0}}$ (if
$t_{k_0} <k$), any $(a_i)_{i=1}^{k_0},\break  (c_i)_{i=1}^{k_0} \in\R
^{k_{0}}$, $(b_i)_{i=1}^{k_0}\in(\R^d)^{k_{0}}$, any $z_{i,j} \in\R
^d$, $\alpha_{i,j}\in\R$, $i=1,\ldots,k_0$, $j=1,\ldots, t_i - t_{i-1}$
(with $t_0= 0$), such that $\| z_{i,j}\| = 1$, $\alpha_{i,j}\geq0$ and
$\sum_{j=1}^{t_i - t_{i-1}}\alpha_{i,j}=1$,
for any $(\gamma_i)_{i=1}^{k - t_{k_0}}$ which belong to $\Gamma
\setminus\{ \gamma_i^0, i=1,\ldots,k_0\}$,
%
%e10 #&#
\begin{equation}
\label{iden:1} \sum_{i=1}^{k - t_{k_0}}
\pi_i g_{\gamma_i} + \sum_{i=1}^{k_0}
\bigl( a_i g_{\gamma_i^0} + b_{i}^{T}D^1g_{\gamma_i^0}
\bigr) + \sum_{i=1}^{k_0}
c_i^2\sum_{j=1}^{t_{i}-t_{i-1}}
\alpha_{i,j} z_{i,j}^{T}D^2
g_{\gamma_i^0} z_{i,j} = 0,
\end{equation}
if and only if
\[
a_i = 0, \qquad b_i = 0,\qquad c_{i}=0
\qquad\forall i=1,\ldots,k_0, \qquad\pi_i=0\qquad
\forall i = 1,\ldots,k-t_{k_0}.
\]
\end{enumerate}

Assumption \textbf{A4} is a weak identifiability condition for
situations when $k_0 < k$.
%, i.e. when the parameter $\theta_0$ is on the boundary of $\Theta_k$.
Notice that \textbf{A4} is the same condition as in Rousseau and
Mengersen \cite{kerrieju:11}, it is satisfied in particular for Poisson mixtures,
location-scale Gaussian mixtures and any mixtures of regular
exponential families.

The following theorem says that the posterior distribution of $L
(\theta )$ concentrates on the true number $k_{0}$ of hidden states.
%An estimator of the order can then be for instance the posterior mode
%of the distribution of $L\left(\theta\right)$, however the whole
%posterior distribution is also of interest.
%%%%%%%%%%%%%%
%
%th2 #&#
\begin{theorem}
\label{theo:order}
Assume that assumptions \textup{\textbf{A0}--\textbf{A2}} and \textup{\textbf{A4}} are
verified. If either of the following two situations holds:
\begin{itemize}
\item Under assumption \textup{\textbf{A3}} (Dirichlet type prior), if $\bar
{\alpha}> k(k+d-1) $ and
\[
\frac{u_n v_n n^{(\bar\alpha- k(k+d-1))/(2\bar\alpha) }}{\log n} \rightarrow+\infty.
\]
\item Under assumption \textup{\textbf{A3bis}} (exponential type prior), if $u_n
v_n n^{1/2 }/(\log n)^{3/2} \rightarrow+\infty$,
\end{itemize}
then
%
%e11 #&#
\begin{equation}
\label{post:k0} \Pb^{\Pi} \bigl[ \theta\dvtx L (\theta )\neq
k_{0} \vert Y_{1:n} \bigr] = \mathrm{o}_{\Pb_{\theta_0}} (1
).
\end{equation}
If $\hat{k}_{n} $ is either the mode or the median of the posterior
distribution of $L(\theta)$, then
%
%e12 #&#
\begin{equation}
\label{pointest:k0} \hat{k}_{n}= k_0 + \mathrm{o}_{\Pb_{\theta_0}}(1).
\end{equation}
\end{theorem}

One of the advantages of using such an estimate of the order of the
HMM, is that we do not need to consider a prior on $k$ and use
reversible-jump methods, see Richardson and Green \cite{richardson:green:1997},
which can
be tricky to implement. In particular, we can consider a two-stage
procedure where $\hat k_{n}$ is computed based on a model with $k$
components where $k$ is a reasonable upper bound\vspace*{1.5pt} on $k_0$ and then,
fixing $k = \hat k_n$ an empirical Bayes procedure is defined on
$(Q_{i,j}, i,j \leq\hat k_n, \gamma_1, \ldots,\gamma_{\hat k_n})$. On
the event $\hat k_n = k_0$, which has probability going to 1 under $\Pb
_{\theta_0}$ the model is regular and using the Bernstein--von Mises
theorem of de~Gunst and Shcherbakova \cite{degunst:shcherbakova:08}, we obtain that with
probability $\Pb_{\theta_0}$ going to 1, the posterior distribution of
$\sqrt{n} ( \theta- \hat\theta_n )$ converges in distribution to the
centered Gaussian with variance $V_0$, the inverse of Fisher
information at parameter $\theta_0$, where $\hat{\theta}_n$ is an
efficient estimator of $\theta_0$ when the order is known to be $k_0$,
and $ \sqrt{n} ( \hat{\theta}_n - \theta_0 )$
converges in distribution to the centered Gaussian with variance $V_0$
under $\Pb_{\theta_0}$.
%[\sqrt{n} ( \theta- \hat\theta_n ) \vert Y_{1:n}] \Rightarrow
%where $\hat\theta_n$ is an efficient estimator conditionally on $\hat
%k = k_0$ and
% $$ \sqrt{n} ( \hat\theta_n - \theta_0 ) \Rightarrow\mathcal N( 0,
%V_0)$$
% under $\Pb_{\theta_0}$. In particular posterior credible regions of
%this empirical Bayes procedure are also asymptotic confidence regions.

The main point in the proof of Theorem~\ref{theo:order} is to prove an
inequality that relates the $L_1$ distance between the $l$-marginals,
to a distance between the parameters of the HMM. Under condition \textbf{A4},
we prove that
there exists a constant $c(\theta_{0})>0$ such that for any small
enough positive $\varepsilon$,
%
%e13 #&#
\begin{eqnarray}
\label{ineqsimple} &&\frac{\llVert  f_{l,\theta}- f_{l,\theta_0} \rrVert _{1}}{c(\theta
_{0})}
\nonumber
\\
&&\quad\geq \sum_{1\leq j \leq k : \forall i, \|\gamma_{j}-\gamma_{i}^{0}\| >
\varepsilon} \Pb_{\theta}
(X_{1}=j ) + \sum_{i=1}^{ k_{0}}
\bigl\llvert \Pb_{\theta} \bigl(X_{1}\in B(i) \bigr)-
\Pb_{\theta_0} (X_{1}=i )\bigr\rrvert
\\
&&\qquad{}+\sum_{i=1}^{ k_{0}} \biggl[\biggl \|\sum
_{j\in B(i)} \Pb_{\theta} (X_{1}=j )
\bigl(\gamma_{j}-\gamma_{i}^{0}\bigr)\biggr\| +
\frac{1}{2} \sum_{j\in B(i)} \Pb_{\theta}
(X_{1}=j )\bigl\llVert \gamma_{j}-\gamma_{i}^{0}
\bigr\rrVert ^{2} \biggr],
\nonumber
\end{eqnarray}
where $B(i)=\{j\dvtx  \|\gamma_{j}-\gamma_{i}^{0}\|\leq\varepsilon\}$. The
above lower bound essentially corresponds to a partition of $\{ 1,\ldots
, k\}$ into $k_0 + 1$ groups, where the first $k_0$ groups correspond
to the components that are close to true distinct components in the
multivariate mixture and the last corresponds to components that are
emptied. The first term on the right-hand side controls the weights of
the components that are emptied (group $k_0 + 1$), the second term
controls the sum of the weights of the components belonging to the
$i$th group, for $i=1,\ldots, k_0$ (components merging with the true
$i$th component), the third term controls the distance between the mean
value over the group $i$ and the true value of the $i$th component in
the true mixture while the last term controls the distance between each
parameter value in group $i$ and the true value of the $i$th component.
A general inequality implying \eqref{ineqsimple}, obtained under a
weaker condition, namely \textbf{A4bis}, holds and is stated and proved
in Section~\ref{subsec:order}.

As we have seen with Theorem~\ref{theo:order}, we can recover the
\textit{true parameter} $\theta_0$ using a two-stage procedure where
first $\hat k_n$ is estimated. However, it is also of interest to
understand better the behaviour of the posterior distribution in the
first stage procedure and see if some behaviour similar to what was
observed in Rousseau and Mengersen \cite{kerrieju:11} holds in the case of HMMs. From Theorem~\ref{theo:gene:finite}, it appears that HMMs present an extra
difficulty due to the fact that, when the order is overestimated, the
neighbourhood of $\theta$'s such that $\Pb_\theta=\Pb_{\theta_0}$
contains parameters leading to non-ergodic HMMs. To have a more refined
understanding of the posterior distribution, we restrict our attention
in Section~\ref{subsec:1et2} to the case where $k=2$ and $k_0 = 1$
which is still nontrivial, see also Gassiat and K\'eribin \cite{CKEG} for the description of
pathological behaviours of the likelihood in such a case.
%%%%%%%%%%%%%%%

%s2.4 #&#
\subsection{Posterior concentration for the parameters: The case $k_{0}=1$ and $k=2$}
\label{subsec:1et2}

In this section, we restrict our attention to the simpler case where
$k_0=1$ and $k=2$. In Theorem~\ref{theo:uncontredeux} below, we prove
that if a Dirichlet type prior is considered on the rows of the
transition matrix with parameters $\alpha_j$'s that are large enough
the posterior distribution concentrates on the configuration where the
two components (states) are merged ($\gamma_1$ and $\gamma_2$ are close
to one another).
When $k=2$, we can parameterize $\theta$ as $\theta=(p,q,\gamma
_{1},\gamma_{2})$, with $0\leq p \leq1$, $0\leq q \leq1$, so that
\[
Q_{\theta}=\pmatrix{ 1-p & p
\cr
q & 1-q }, \qquad\mu_{\theta}=
\biggl(\frac{q}{p+q}, \frac{p}{p+q} \biggr)
\]
when $p\neq0$ or $q\neq0$. If $p=0$ and $q=0$, set $\mu_{\theta}=
(\frac{1}{2}, \frac{1}{2} )$, for instance. Also, we may take
\[
\rho_{\theta}-1= (p+q ) \wedge \bigl(2-(p+q) \bigr).
\]
When $k_{0}=1$, the observations are i.i.d. with distribution $g_{\gamma
^{0}}\,\mathrm{d}\nu$, so that one may take
$\theta_{0}=(p,1-p,\gamma^{0},\gamma^{0})$ for any $0<p<1$, or
$\theta_{0}=(0,q,\gamma^{0},\gamma)$ for any $0<q\leq1$ and any $\gamma
$, or
$\theta_{0}=(p,0,\gamma,\gamma^{0})$ for any $0<p\leq1$ and any $\gamma$.
Also, for any $x\in{\mathcal X}$, $\Pb_{\theta_{0},x}=\Pb_{\theta
_{0}}$ and
\[
\ell_n(\theta,x) - \ell_n(\theta_0,x_{0})=
\ell_n(\theta,x) - \ell _n(\theta_0,x).
\]
We take independent Beta priors on $(p,q)$:
\[
\Pi_{2} (\mathrm{d}p,\mathrm{d}q)= C_{\alpha,\beta} p^{\alpha-1}(1-p)^{\beta-1}q^{\alpha
-1}(1-q)^{\beta-1}
\one_{0<p<1}\one_{0<q<1}\,\mathrm{d}p \,\mathrm{d}q,
\]
thus satisfying \textbf{A3}.
%Let $\delta_{n}=\frac{M_{n}}{\sqrt{n}}$ with
%We will sethanks to Theorem~\ref{theo:gene:finite},
%Let $\epsilon_{n} $ by any sequence converging to 0 and consider the
%set
%$$
%A_{n}=\left\{\frac{p}{p+q} \leq\epsilon_{n} {\text{or }}
%$$
Then the following holds.
%
%th3 #&#
\begin{theorem}
\label{theo:uncontredeux}
Assume that assumptions \textup{\textbf{A0}--\textbf{A2}} together with
assumption \textup{\textbf{A4}} are verified and consider the prior described
above with $\omega(\cdot)$ verifying \textup{\textbf{A3}}. Assume moreover that
for all $x$, $\gamma\mapsto g_{\gamma}(x)$ is four times continuously
differentiable on $\Gamma$, and that for any $\gamma\in\Gamma$ there
exists $\epsilon>0$ such that for any $i \leq4$,
%
%e14 #&#
\begin{equation}
\label{hypocor2} \int\sup_{\gamma'\in B_d(\gamma, \epsilon)}\biggl\| \frac{D^{i}_{\gamma}
g_{\gamma'}}{g_{\gamma'}} (y )
\biggr\|^{4}g_{\gamma} (y )\nu(\mathrm{d}y)< +\infty.
\end{equation}
%
%where $D^{i}_{\gamma}g_{\gamma'}$ denotes the $i$-th differential
%operator (with respect to $\gamma$) of $g$ at point $\gamma'$,
%for any sequence $(\epsilon_{n})_{n\geq1}$ decreasing slowly to $0$,
Then, as soon as $\alpha>3d/4$ and $\beta>3d/4$, for any sequence
$\epsilon_{n}$ tending to $0$,
\[
\Pb^{\Pi} \biggl(\frac{p}{p+q} \leq\epsilon_{n} \mbox{
or } \frac
{q}{p+q} \leq\epsilon_{n}\big\vert Y_{1:n}
\biggr)= \mathrm{o}_{\Pb_{\theta_0}}(1),
\]
and for any sequence $M_n$ going to infinity,
\[
\Pb^{\Pi} \bigl( \|\gamma_1 - \gamma_0\| + \|
\gamma_2 - \gamma_0\| \leq M_n
n^{-1/4} \big\vert Y_{1:n} \bigr) = 1 + \mathrm{o}_{\Pb_{\theta_0}}(1).
\]
\end{theorem}

Theorem~\ref{theo:uncontredeux} says that the extra component cannot be
emptied at rate $\epsilon_{n} $, where the sequence $\epsilon_n$ can be
chosen to converge to 0 as slowly as we want, so that asymptotically,
under the posterior distribution neither $p/(p+q)$ nor $q/(p+p)$ are
small, and the posterior distribution concentrates on the configuration
where the components merge, with the emission parameters merging at
rate $n^{-1/4}$. Similarly in Rousseau and Mengersen \cite{kerrieju:11}
the authors obtain
that, for independent variables, under a Dirichlet $\mathcal{D}(\alpha
_1,\ldots, \alpha_k)$ prior on the weights of the mixture and if $\min
\alpha_i > d/2$, the posterior distribution concentrates on
configurations which do not empty the extra-components but merge them
to true components. The threshold here is $ 3d/2$ instead of $d/2$.
This is due to the fact that there are more parameters involved in a
HMM model associated to $k$ states than in a $k$-components mixture
model. No result is obtained here in the case where the $\alpha_i$'s
are small. This is due to the existence of non ergodic $\Pb_\theta$ in
the vicinity of $\Pb_{\theta_0}$ that are not penalized by the prior in
such cases. Our conclusion is thus to favour large values of the $\alpha_i$'s.

%%%%%%%%%%%%%

%s3 #&#
\section{A general theorem}
\label{sec:gene}

In this section, we present a general theorem which is used to prove
Theorem~\ref{theo:gene:finite} but which can be of interest in more
general HMMs.
We assume here that the unobserved Markov chain $(X_i)_{i\geq1}$ lives
in a Polish space $\mathcal X$ and the observations $(Y_i)_{i\geq1}$
are conditionally independent given $(X_i)_{i\geq1}$ and live in a
Polish space $\mathcal Y$. $\mathcal X$, $\mathcal Y$ are endowed with
their Borel $\sigma$-fields.
We denote by $\theta\in\Theta$, where $\Theta$ is a subset of an
Euclidean space, the parameter describing the distribution of the HMM,
so that $Q_\theta$, $\theta\in\Theta$ is the Markov kernel of
$(X_i)_{i\geq1}$ and the conditional distribution of $Y_i $ given
$X_i$ has density with respect to some given measure $\nu$ on $\mathcal
Y$ denoted by $g_\theta(y|x)$, $x\in\mathcal X$, $\theta\in\Theta$.
%With an abuse of notations we also denote $\nu$ the product measure $
We assume that the Markov kernels $Q_\theta$ admit a (not necessarily
unique) stationary distribution $\mu_\theta$, for each $\theta\in
\Theta$.
We still write $\Pb_{\theta}$ for the probability distribution of the
stationary HMM $(X_{j},Y_{j})_{j\geq1}$ with parameter $\theta$. That
is, for any integer $n$, any set $A$ in the Borel $\sigma$-field of
${\mathcal X}^{n}\times{\mathcal Y}^{n}$:
%
%e15 #&#
\begin{eqnarray}
\label{defhmmstat} &&\Pb_{\theta} \bigl((X_{1},
\ldots,X_{n},Y_{1},\ldots,Y_{n})\in A \bigr)
\nonumber
\\[-8pt]
\\[-8pt]
&&\quad=\int_{A}\mu_{\theta}(\mathrm{d}x_{1})
\prod_{i=1}^{n-1}Q_{\theta}
(x_{i},\mathrm{d}x_{i+1} )\prod
_{i=1}^{n}g_{\theta} (y_{i}\vert
x_{i} )\nu(\mathrm{d}y_{1})\cdots\nu(\mathrm{d}y_{n}).
\nonumber
\end{eqnarray}
Thus for any integer $n$, under $\Pb_{\theta}$, $Y_{1:n}=(Y_{1},\ldots
,Y_{n})$ has a probability density with respect to $\nu(\mathrm{d}y_{1})\cdots\nu
(\mathrm{d}y_{n})$ equal to
%
%e16 #&#
\begin{equation}
\label{margeq2} f_{n,\theta}(y_1,\ldots, y_n) =
\int_{{\mathcal X}^{n}}\mu_{\theta}(\mathrm{d}x_{1})
\prod_{i=1}^{n-1}Q_{\theta}
(x_{i},\mathrm{d}x_{i+1} )\prod
_{i=1}^{n}g_{\theta} (y_{i}\vert
x_{i} ).
\end{equation}
We denote by $\Pi_{\Theta}$ the prior distribution on $\Theta$ and by
$\pi_{\mathcal X}$ the prior probability on the unobserved initial
state, which might be different from the stationary distribution $\mu
_\theta$.
We set $\Pi=\Pi_{\Theta}\otimes\pi_{\mathcal X}$. Similarly to before,
denote $\ell_n(\theta,x)$ the log-likelihood starting from $x$, for all
$x \in\mathcal X$.
%, %which is given by
%$$
%$$
%Similarly, is denoted
%and by $\ell_{n}\left(\theta,\pi_0\right)$ the log-likelihood starting
%from a distribution $\pi_0$ on $\mathcal X$.
%$$\ell_{n}\left(\theta,\pi_0\right) = \log\left[ \int_{\mathcal X} e^{

%The posterior distribution can then be written as
% \end{equation}
%for any Borel set $A\subset\Theta$.\\
%We shall also use the notation $\Pb_{\theta,x}$ for the probability
%distribution of the HMM starting from $x$, that is, for any integer
%$n$, any measurable set $A$ in the Borel $\sigma$-field of ${\mathcal
%X}^{n}\times{\mathcal Y}^{n}$:
%so that for any $\theta\in\Theta$,
%$$
%$$
%We denote by $E_{\theta}$ the expectation under $\Pb_{\theta}$ and by
%$E_{\theta,x}$ the expectation under $\Pb_{\theta,x}$.\\

We assume that we are given a stationary HMM $(X_{j},Y_{j})_{j\geq1}$
with distribution $\Pb_{\theta_{0}}$ for some $\theta_{0}\in\Theta$.

For any $\theta\in\Theta$,
%since the total variation norm between probability measures is bounded
%by $2$,
it is possible to define real numbers $\rho_{\theta}\geq1$
and
$0<\R_{\theta}\leq2$ such that, for any integer $m$, any $x\in
{\mathcal X}$
%
%e17 #&#
\begin{equation}
\label{geo:ergo} \bigl \| Q_\theta^m (x,\cdot) - \mu_\theta
\bigr \|_{\mathrm{TV}} \leq R_\theta\rho_{\theta}^{-m},
\end{equation}
where $\|\cdot\|_{\mathrm{TV}}$ is the total variation norm. If it is possible
to set $\rho_{\theta}>1$, the Markov chain $(X_{n})_{n\geq1}$ is
uniformly ergodic and $\mu_{\theta}$ is its unique stationary
distribution. The following theorem provides a posterior concentration
result in a general HMM setting, be it parametric or nonparametric and
is an adaptation of Ghosal and van~der Vaart \cite{ghosal:vaart:2007} to the setup of HMMs. We
present the assumptions needed to derive the posterior concentration rate.
\begin{enumerate}[\textbf{C3bis}]
\item[\textbf{C1}] There exists $A>0$ such that for any
$(x_{0},x_{1})\in\mathcal{X}^{2}$, $\Pb_{\theta_0}$ almost surely,
$\forall n \in\N$,
$| \ell_n(\theta_0, x_{0})-\ell_n(\theta_0, x_{1})| \leq A$,
and
there exist $S_n \subset\Theta\times\mathcal X$, $C_n >0$ and
$\tilde\epsilon_n>0$ a sequence going to 0 with $n \tilde\epsilon_n^2
\rightarrow+\infty$
such that
\[
\sup_{(\theta,x) \in S_{n}} \Pb_{\theta_0} \bigl[ \ell_n(
\theta,x) - \ell_n(\theta_0, x_{0}) \leq- n
\tilde\epsilon_n^2 \bigr] = \mathrm{o}(1),\qquad
\Pi[S_n] \gtrsim \mathrm{e}^{- C_n n \tilde\epsilon_n^2 }.
\]

\item[\textbf{C2}] There exists a sequence $(\mathcal F_n)_{n\geq1}$
of subsets of $\Theta$
\[
\Pi_{\Theta}\bigl(\mathcal F_{n}^c \bigr) =
\mathrm{o}\bigl(\mathrm{e}^{-n\tilde\epsilon_n^2 (1 + C_n)}\bigr).
\]

\item[\textbf{C3}] There exists a sequence $\epsilon_n \geq\tilde
\epsilon_n$ going to 0, such that $(n\tilde\epsilon_n^2 (1 +
C_n))/(n\epsilon_{n}^{2})$ goes to $0$ and
\[
N \biggl( \frac{\epsilon_n }{12}, \mathcal F_n, d_{l}(\cdot,
\cdot) \biggr)\leq \mathrm{e}^{ \trup{(n \epsilon_n^2 ( \rho_{\theta_0}- 1)^2 )}{( 16 l (
2R_{\theta_0} + \rho_{\theta_0}-1)^2 )}},
\]
where $N(\delta, \mathcal F_n, d_{l}(\cdot,\cdot) )$ is the smallest
number of $\theta_j \in\mathcal F_n$ such that for all $\theta\in
\mathcal F_n$ there exists a $\theta_j$ with $d_{l}(\theta_j, \theta
)\leq\delta$.

Here $d_{l}(\theta,\theta_j) = \|f_{l,\theta} - f_{l,\theta_j}\|
_1:=\int_{\mathcal Y^l}|f_{l,\theta} - f_{l,\theta_j}|(y)\,\mathrm{d}\nu^{\otimes l}(y)$.

\item[\textbf{C3bis}] There exists a sequence $\epsilon_n \geq\tilde
\epsilon_n$ going to 0 such that
\[
\sum_{m \geq1} \frac{\Pi_{\Theta}  (A_{n,m} (\epsilon_{n}
) )}{ \Pi ( S_n  ) }\mathrm{e}^{ - \trup{(nm^2\epsilon_n^2 )}{( 32 l
)} }
= \mathrm{o}\bigl(\mathrm{e}^{- n\tilde\epsilon_n^2}\bigr)
\]
and
\[
N \biggl(\frac{m\epsilon_{n}}{12}, A_{n,m} (\epsilon_{n} ),
d_{l}(\cdot,\cdot) \biggr) \leq \mathrm{e}^{ \trup{(n m^2 \epsilon_n^2 ( \rho
_{\theta_0}- 1)^2 )}{(16 l ( 2R_{\theta_0} + \rho_{\theta_0}-1)^2)}},
\]
where
\[
A_{n,m} (\epsilon )= \mathcal F_n \cap \biggl\{\theta\dvtx
m\epsilon\leq\|f_{l,\theta} - f_{l,\theta_0}\|_1
\frac{\rho_\theta-
1}{2R_{\theta}+\rho_\theta- 1} \leq(m+1)\epsilon \biggr\}.
\]
\end{enumerate}

%th4 #&#
\begin{theorem}\label{theo:gene}
Assume that $\rho_{\theta_0}>1$ and that assumptions \textup{\textbf{C1}--\textbf{C2}} are
satisfied, together with either assumption \textup{\textbf{C3}} or
\textup{\textbf{C3bis}}. Then
\[
\Pb^\Pi \biggl[ \theta\dvtx \|f_{l,\theta} - f_{l,\theta_0}
\|_1 \frac{\rho
_\theta- 1}{2R_{\theta}+\rho_\theta- 1} \geq \epsilon_n \Big\vert
Y_{1:n} \biggr] = \mathrm{o}_{\Pb_{\theta_0}} (1 ).
\]
\end{theorem}

Theorem~\ref{theo:gene} gives the posterior concentration rate of $ \|
f_{l,\theta} - f_{l,\theta_0}\|_1 $ up to the parameter $\frac{\rho
_\theta- 1}{2R_{\theta}+\rho_\theta- 1}$.
In Ghosal and van~der Vaart \cite{ghosal:vaart:2007}, for models of non independent variables,
the authors consider a parameter space where the mixing coefficient
term (for us $\rho_\theta- 1$) is uniformly bounded from below by a
positive constant over $\Theta$
(see their assumption (4.1) for the application to Markov chains or
their assumption on $\mathcal F$ in Theorem~7 for the application to
Gaussian time series),
%and where $R_\theta$ is uniformly bounded
or equivalently they consider a prior whose support in $\Theta$ is
included in a set where $\frac{\rho_\theta- 1}{2R_{\theta}+\rho_\theta
- 1}$ is uniformly bounded from below, so that their posterior
concentration rate is directly expressed in terms of $ \|f_{l,\theta} -
f_{l,\theta_0}\|_1$. Since we do not restrict ourselves to such
frameworks the penalty term $\rho_\theta-1$ is incorporated in our
result. However Theorem~\ref{theo:gene}, is proved along the same lines
as Theorem~1 of Ghosal and van~der Vaart \cite{ghosal:vaart:2007}.

%In this paper we apply Theorem~\ref{theo:gene} to finite state space
%HMMs where the use of truncated priors over sets where $\frac{\rho_
%below is awkward and for which this term may be seen as a penalty term
%and plays a role in the study of asymptotic posterior concentration.
%In particular in the case of over-fitted HMMs with finite state space,
%i.e. when $\theta_0$ corresponds to a HMM associated with $k_0$ states
%while the model considers HMMs with $k>k_0$ states, the parameter set
%has to contain all possible transition matrices and in any
%neighbourhood $ \{ \theta:\| f_{l,\theta} - f_{l,\theta_0}\|_1\leq
%there exist parameters $\theta$ such that $\rho(\theta)=1$. Thus, one
%has to allow $\rho_{\theta}$ to be arbitrarily close to $1$.
%The prior, however, acting as soft thresholding, leads to the
%concentration of the posterior distribution around $f_{\theta_0}$, in
%terms of $ \|f_{l,\theta} - f_{l,\theta_0}\|_1 $ alone, at a rate
%slower than $(\log n/n)^{1/2}$, see Theorem~\ref{theo:gene:finite}.

The assumption $\rho_{\theta_0}>1$ implies that the hidden Markov chain
$X$ is uniformly ergodic. Assumptions \textbf{C1}--\textbf{C2} and either
\textbf{C3} or \textbf{C3bis}
are similar in spirit to those considered in general theorems on
posterior consistency or posterior convergence rates, see, for
instance, Ghosh and Ramamoorthi \cite{ghosh:ramamo:2003} and Ghosal and
van~der Vaart \cite{ghosal:vaart:2007}.
Assumption \textbf{C3bis} is often used to eliminate some extra $\log
n$ term which typically appear in nonparametric posterior concentration
rates and is used in particular in the proof of Theorem~\ref{theo:uncontredeux}.

%%%%%%%%%%%%%%%%%%%

%s4 #&#
\section{Proofs}
\label{sec:proofs}

%s4.1 #&#
\subsection{Proof of Theorem \texorpdfstring{\protect\ref{theo:gene:finite}}{1}}
\label{subsec:gene:finite}

The proof consists in showing that the assumptions of
Theorem~\ref{theo:gene} are satisfied.

Following
the proof of Lemma~2 of Douc \textit{et al.} \cite{douc:moulines:ryden:04} we find that,
since $\rho_{\theta_0}>1$, for any $x_{0}\in\mathcal{X}$,
\[
\bigl | \ell_{n} (\theta_{0},x_{0} )-
\ell_{n} (\theta _{0},x_{1} )\bigr | \leq2 \biggl(
\frac{\rho_{\theta_{0}}}{\rho_{\theta_{0}}-1} \biggr)^{2}
\]
so that setting $A=2 (\frac{\rho_{\theta_{0}}}{\rho_{\theta
_{0}}-1} )^{2}$ the first point of \textbf{C1} holds.

We shall verify assumption \textbf{C1} with $\tilde\epsilon_n= M_{n} /
\sqrt{n}$ for some $M_{n}$ tending slowly enough to infinity and that
will be chosen later.
Note that
the assumption \textbf{A0} and the construction (\ref{lem:m1}) allow to
define a $\tilde{\theta}_{0}\in\Theta_{k}$ such that, writing $\tilde
{\theta}_{0} = (\tilde Q^0, \tilde\gamma_1^0, \ldots, \tilde\gamma
_k^0)$ with $\tilde Q^0 = (\tilde q_{i,j}^0, i,j\leq k)$,
%$\rho_Q> 1$.
%Assumption \textbf{A0} and the construction (\ref{lem:m1}) allow to
%define a $\tilde{\theta}_{0}\in\Theta_{k}$ such that \textbf{A0} holds
%with $D=k(k-1)+kd$. We write $\theta_0 = (Q^0, \gamma_1^0, \ldots,
if $V$ is a bounded subset of $\{ \theta= (Q, \gamma_1,\ldots, \gamma
_k); |q_{i,j}-\tilde q_{i,j}^0| \leq\tilde\epsilon_n\}$, then
%
%e18 #&#
\begin{equation}
\label{P0} \inf_{\theta\in V} \rho_\theta> 1,
\end{equation}
for large enough $n$, and
\[
\sup_{\theta\in V}\sup_{x,x_{0}\in{\mathcal X}} \bigl\llvert \ell
_{n} (\theta,x )-\ell_{n} (\theta,x_{0} ) \bigr
\rrvert \leq2 \sup_{\theta\in V} \biggl(\frac{\rho_{\theta}}{\rho
_{\theta}-1}
\biggr)^{2}.
\]
Following the proof of Lemma~2 of Douc \textit{et al.} \cite{douc:moulines:ryden:04} gives
that, if \textbf{A0} and \textbf{A1} hold, for all $\theta\in V$ $\Pb
_{\theta_0}$-a.s.,
%
%e19 #&#
\begin{eqnarray}
\label{lik2} %
 \ell_{n} (\theta,x_{0} )-
\ell_{n} (\theta_{0},x_{0} )
&=& (\theta-\theta_0 )^{T}\nabla_{\theta}
\ell_{n} (\theta _{0},x_{0} )
\nonumber
\\[-8pt]
\\[-8pt]
&&{} +\int_{0}^{1} (\theta-
\theta_0 )^{T}D^{2}_{\theta}\ell
_{n} \bigl(\theta_{0}+u(\theta-\theta_{0}),x_{0}
\bigr) (\theta -\theta_0 ) (1-u)\,\mathrm{d}u.
\nonumber
\end{eqnarray}
Following Theorem~2 in Douc \textit{et al.} \cite{douc:moulines:ryden:04},
$n^{-1/2}\nabla_\theta\ell_n(\theta_0, x)$ converges in distribution
under $\Pb_{\theta_0}$ to ${\mathcal N}(0, V_0)$
for some positive definite matrix $V_0$, and following Theorem~3 in Douc \textit{et al.}
\cite{douc:moulines:ryden:04}, we get that
$ \sup_{\theta\in V} n^{-1} D_\theta^2 \ell_n(\theta, x_{0}) $
converges $\Pb_{\theta_0}$ a.s. to $V_0$.
Thus, we may set:
\[
S_{n}= \bigl\{ \theta\in V; \bigl \| \gamma_j -
\gamma_j^0\bigr \| \leq1/ \sqrt{n} \ \forall j \leq k \bigr\}
\times{\mathcal X}
\]
so that
%
%e20 #&#
\begin{equation}
\label{denom} \sup_{(\theta, x)\in S_n} \Pb_{\theta_0} \bigl[
\ell_n(\theta,x) - \ell _n(\theta_0,
x_0) < -M_n \bigr] = \mathrm{o}(1).
\end{equation}
Moreover, letting $D = k(k-1+d)$, we have $\Pi\otimes\Pi_{\mathcal
{X}}(S_n) \gtrsim n^{-D/2}$ and \textbf{C1} is then satisfied setting
$C_n = D\log n/(2M_{n}^{2})$.

% To prove that \textbf{C2} and \textbf{C3} hold, recall that if $
Let now $v_{n}=n^{-D/(2\min_{1\leq i \leq k}\alpha_{i})}/\sqrt{\log n}$
and $u_{n}=n^{-D/(2\sum_{1\leq i \leq k}\alpha_{i})}/\sqrt{\log n}$,
and define
%Let us now define, for sequences $u_{n}$ and $v_{n}$
\begin{eqnarray*}
 {\mathcal F}_{n}&=& \Biggl\{\theta= (q_{ij}, 1\leq i \leq k, 1
\leq j \leq k-1 ; \gamma_1,\ldots,\gamma_k ) \dvtx
q_{ij}\geq v_{n}, 1\leq i \leq k, 1\leq j \leq k,
\\
&&\phantom{\Biggl\{}  \sum
_{j=1}^{k}\inf_{1\leq i \leq k}q_{ij}
\geq u_{n}, \| \gamma _i\| \leq{n}^{b}, 1\leq
i \leq k \Biggr\}.
\end{eqnarray*}
Now, if $\theta\in\mathcal{F}_{n}^c$, then there exist $1\leq i,j \leq
k$ such that $q_{ij}\leq v_{n}$, or $\sum_{j=1}^{k}\inf_{1\leq i \leq
k}q_{ij} \leq u_{n}$, or
there exists $ 1\leq i \leq k$ such that $\| \gamma_i\| \geq{n}^{b}$. Using
\textbf{A3} we easily obtain that for fixed $i$ and $j$, $\Pi(\{\theta
\dvtx q_{ij}\leq v_{n}\})=\mathrm{O}(v_{n}^{\alpha_{j}})$ and $\Pi(\{\theta\dvtx \| \gamma
_i\| \geq{n}^{b}\})=\mathrm{o}(n^{-D/2}) $.
Also, if $\sum_{j=1}^{k}\inf_{1\leq i \leq k}q_{ij} \leq u_{n}$, then
there exists a function $i(\cdot)$ from $\{1,\ldots,k\}$ to $\{1,\ldots
,k\}$ whose image set has cardinality at least $2$ such that $\sum_{j=1}^{k} q_{i(j)j} \leq u_{n}$.
This gives, using \textbf{A3}, $\Pi(\{\theta\dvtx  \sum_{j=1}^{k}\inf_{1\leq i \leq k}q_{ij} \leq u_{n}\})=\mathrm{O}(u_{n}^{\sum_{1\leq i \leq
k}\alpha_{i}})$.
Thus,
\[
\Pi\bigl(\mathcal F_{n}^c \bigr) = \mathrm{O}
\bigl(v_{n}^{\min_{1\leq i \leq k}\alpha
_{i}}+u_{n}^{\sum_{1\leq i \leq k}\alpha_{i}}\bigr)+
\mathrm{o} \bigl(n^{-D/2} \bigr).
\]
We may now choose $M_{n}$ tending to infinity slowly enough so that
$v_{n}^{\min_{1\leq i \leq k}\alpha_{i}}+u_{n}^{\sum_{1\leq i \leq
k}\alpha_{i}}=\mathrm{o}(\mathrm{e}^{-M_{n}}n^{-D/2})$ and $\Pi(\mathcal F_{n}^c
)=\mathrm{o}(\mathrm{e}^{-M_{n}}n^{-D/2})$. Then, \textbf{C2} holds.

Now, using the definition of $f_{l,\theta}$, we obtain that
\[
\llVert f_{l,\theta_1}-f_{l,\theta_2} \rrVert _{1} \leq\sum
_{j=1}^k |\mu_{\theta_1}-
\mu_{\theta_2}| + l \sum_{i,j=1}^k
\bigl |Q_{i,j}^1-Q_{i,j}^2\bigr |+ l \max
_{j \leq k} \|g_{\gamma_j^1}- g_{\gamma_j^2}\|_1
\]
so that using Lemma~\ref{lem:der} below, \textbf{A1} and \textbf{A2} we
get that
for some constant $B$, $\forall(\theta_{1},\theta_{2})\in{\mathcal F}_{n}^{2}$
\[
\llVert f_{l,\theta_1}-f_{l,\theta_2} \rrVert _{1}\leq B
\biggl(\frac{1}{v_{n}^{2c}}+n^{a} \biggr) \llVert \theta_{1}-
\theta _{2}\rrVert .
\]
Thus for some other constant $\tilde B$,
\[
N\bigl(\delta, \mathcal F_n, d(\cdot,\cdot) \bigr)\leq \biggl[
\frac{\tilde
B}{\delta} \biggl(\frac{1}{v_{n}^{2c}}+n^{a} \biggr) \biggr]
^{k(k-1)+kd}
\]
and \textbf{C3} holds when setting $\epsilon_{n}=K\sqrt{\frac{\log
n}{n}}$ with $K$ large enough.

We have proved that under assumptions \textbf{A0}, \textbf{A1},
\textbf{A2}, \textbf{A3}, Theorem~\ref{theo:gene} applies with $\epsilon
_{n}=K\sqrt{\frac{\log n}{n}}$ so that
\[
\Pb^{\Pi} \biggl[  \|f_{l,\theta} - f_{l,\theta_0}
\|_1(\rho_\theta- 1) \geq K \sqrt{\frac{\log n}{n}}\Big\rrvert
Y_{1:n} \biggr] =\mathrm{o}_{\Pb
_{\theta_0}} (1 )
\]
and the first part of Theorem~\ref{theo:gene:finite} is proved. Now
\begin{eqnarray*}
\mathrm{o}_{\Pb_{\theta_0}} (1 ) &=& \Pb^{\Pi} \biggl[
\|f_{l,\theta
} - f_{l,\theta_0}\|_1(\rho_\theta- 1)
\geq K \sqrt{\frac{\log n}{n}} \Big\rrvert Y_{1:n} \biggr]
\\
&=& \Pb^{\Pi} \biggl[  \theta\in{\mathcal F}_{n}
\mbox{ and }\| f_{l,\theta} - f_{l,\theta_0}\|_1 (
\rho_\theta- 1) \geq K \sqrt{\frac
{\log n}{n}}\Big\rrvert Y_{1:n}
\biggr] +\mathrm{o}_{\Pb_{\theta_0}} (1 ).
\end{eqnarray*}
Since $\rho_\theta- 1 \geq\sum_{j=1}^{k}\min_{1\leq i \leq k}q_{ij}$,
for all $ \theta\in{\mathcal F}_{n}$, $\rho_\theta- 1\geq u_{n}$,
\[
 \Pb^{\Pi} \biggl[ \|f_{l,\theta} - f_{l,\theta_0}
\|_1 (\rho_\theta- 1) \geq K \sqrt{\frac{\log n}{n}} \Big\vert
Y_{1:n} \biggr]
 \geq \Pb^{\Pi} \biggl[ \|f_{l,\theta} - f_{l,\theta_0}
\|_1 \geq2K \frac
{1}{u_{n}}\sqrt{\frac{\log n}{n}} \Big\vert
Y_{1:n} \biggr],
\]
and the theorem follows when \textbf{A3} holds. If now \textbf{A3bis}
holds instead of \textbf{A3}, one gets, taking $u_{n}=v_{n} = h/ \log
n$, with $h > 2C/(k+d-1)$
\[
\Pi\bigl(\mathcal F_{n}^c \bigr) = \mathrm{O}
\bigl(v_{n}\exp(-C/v_{n})\bigr)+ \mathrm{o}
\bigl(n^{-D/2} \bigr) = \mathrm{o}\bigl(\mathrm{e}^{-M_n} n^{-D/2}
\bigr)
\]
by choosing $M_n$ increasing to infinity slowly enough so that
\textbf{C2} and \textbf{C3} hold. The end of the proof follows similarly as before.

To finish the proof of Theorem~\ref{theo:gene:finite},
we need to prove the following lemma.
%
%le1 #&#
\begin{lemma}
\label{lem:der}
The function $\theta\mapsto\mu_{\theta}$ is continuously
differentiable in $( \Delta_{k}^{0})^{k}\times\Gamma^{k}$ and there
exists an integer $c>0$ and a constant $C>0$ such that for any $1 \leq
i \leq k$, $1 \leq j \leq k-1$, any $m=1,\ldots,k$,
\[
\biggl\llvert \frac{\partial\mu_{\theta} (m )}{\partial q_{ij}} \biggr\rrvert \leq\frac{C}{(\inf_{i'\neq j'}q_{i'j'} )^{2c}}.
\]
One may take $c=k-1$.
\end{lemma}

Let $\theta = (q_{ij}, 1\leq i \leq k, 1\leq j \leq k-1 ; \gamma
_1,\ldots,\gamma_k ) $ be such that $ (q_{ij}, 1\leq i \leq k, 1\leq
j \leq k-1)\in\Delta_{0}^{k}$,
$Q_{\theta}= (q_{ij}, 1\leq i \leq k, 1\leq j \leq k)$ is a $k\times k$
stochastic matrix with positive entries, and $\mu_{\theta}$ is uniquely
defined by the equation
\[
\mu_{\theta}^{T} Q_{\theta} = \mu_{\theta}^{T}
\]
if $\mu_{\theta}$ is the vector $(\mu_{\theta}(m))_{1\leq m\leq k}$.
This equation is solved by linear algebra as
%
%e21 #&#
\begin{eqnarray}
\label{eq:der0} \mu_{\theta} (m )&=&\frac{P_{m}(q_{ij}, 1\leq i \leq k, 1\leq j
\leq k-1)}{R(q_{ij}, 1\leq i \leq k, 1\leq j \leq k-1)},\qquad m=1,\ldots
,k-1,
\nonumber
\\[-8pt]
\\[-8pt]
\mu_{\theta} (k )&=&1-\sum_{m=1}^{k-1}
\mu_{\theta} (m ),
\nonumber
\end{eqnarray}
where $P_{m}$, $l=1,\ldots,k-1$ and $R$ are polynomials where the
coefficients are integers (bounded by $k$) and the monomials are all of
degree $k-1$, each variable
$q_{ij}$, $1\leq i \leq k$, $1\leq j \leq k-1$ appearing with power $0$ or
$1$. Now, since the equation has a unique solution as soon as $
(q_{ij}, 1\leq i \leq k, 1\leq j \leq k-1)\in\Delta_{0}^{k}$, then $R$
is never $0$ on $\Delta_{0}^{k}$, so it may be $0$ only at the
boundary. Thus, as a fraction of polynomials with nonzero denominator,
$\theta\mapsto\mu_{\theta}$ is infinitely differentiable in $( \Delta
_{k}^{0})^{k}\times\Gamma^{k}$, and the derivative has components all
of form
\[
\frac{P(q_{ij}, 1\leq i \leq k, 1\leq j \leq k-1)}{R(q_{ij}, 1\leq i
\leq k, 1\leq j \leq k-1)^{2}},
\]
where again $P$ is a polynomial where the coefficients are integers
(bounded by $2k$) and the monomials are all of degree $k-1$, each variable
$q_{ij}$, $1\leq i \leq k$, $1\leq j \leq k-1$ appearing with power $0$ or
$1$. Thus, since all $q_{ij}$'s are bounded by $1$ there exists a
constant $C$ such that for all $m=1,\ldots,k$, $i=1,\ldots,k$,
$j=1,\ldots,k-1$,
%
%e22 #&#
\begin{equation}
\label{eq:der1} \biggl\llvert \frac{\partial\mu_{\theta}(m)}{\partial q_{ij}} \biggr\rrvert \leq
\frac{C}{R(q_{ij}, 1\leq i \leq k, 1\leq j \leq k-1)^{2}}.
\end{equation}
We shall now prove that
%
%e23 #&#
\begin{equation}
\label{eq:der2} R (q_{ij}, 1\leq i \leq k, 1\leq j \leq k-1 ) \geq
\Bigl(\inf_{1\leq i \leq k, 1\leq j \leq k, i\neq j}q_{ij}\Bigr)^{k-1},
\end{equation}
which combined with (\ref{eq:der1}) and
(\ref{eq:der2}) implies Lemma~\ref{lem:der}.
Note that we can express $R$ as a polynomial function of $Q = q_{ij}$,
$1\leq i \leq k, 1\leq j \leq k, i\neq j$. Indeed, $\mu:=(\mu_{\theta
}(i))_{1\leq i \leq k-1}$ is solution of
\[
\mu^{T} \cdot M = V^{T},
\]
where $V$ is the $(k-1)$-dimensional vector $(q_{kj})_{1\leq j \leq
k-1}$, and $M$ is the $(k-1)\times(k-1)$-matrix with components
$M_{i,j}=q_{kj}-q_{ij}+\one_{i=j}$. Since $R$ is the determinant of
$M$, this leads to, for any $k\geq2$:
%
%e24 #&#
\begin{equation}
\label{eqdet} R=\sum_{\sigma\in{\mathcal S}_{k-1}} \varepsilon ( \sigma )
\prod_{1\leq i \leq k-1, \sigma(i)=i} \biggl(q_{ki}+\sum
_{1\leq j \leq
k-1, j\neq i}q_{ij} \biggr) \prod
_{1\leq i \leq k-1, \sigma(i)\neq i} (q_{ki}-q_{\sigma
(i)i} ),
\end{equation}
where for any integer $n$, ${\mathcal S}_{n}$ is the set of
permutations of $\{1,\ldots,n\}$, and for each permutation $\sigma$, $
\varepsilon ( \sigma )$ is its signature. Thus,
$R$ is a polynomial in the components of $Q$ where each monomial has
integer coefficient and has $k-1$ different factors. The possible
monomials are of form
\[
\beta\prod_{i\in A} q_{ki} \prod
_{i\in B} q_{i j(i)},
\]
where $(A,B)$ is a partition of $\{1,\ldots,k-1\}$, and for all $i\in
B$, $j(i)\in\{1,\ldots,k-1\}$ and $j(i)\neq i$. In case $B=\emptyset$,
the coefficient $\beta$ of the monomial is $\sum_{\sigma\in{\mathcal
S}_{k-1}} \varepsilon ( \sigma )=0$, so that we only consider
partitions
such that $B\neq\emptyset$.  Fix such a monomial with non-null
coefficient, let $(A,B)$ be the associated partition. Let $Q$ be such
that, for all $i\in A$, $q_{ki}>0$, for all $i\notin A$, $q_{ki}=0$ and
$q_{kk}>0$ (used to handle the case $A=\emptyset$). Fix also
$q_{ij(i)}=1$ for all $i\in B$. Then, if $(A',B')$ is another partition
of $\{1,\ldots,k-1\}$ with $B'\neq\emptyset$,
the monomial $ \prod_{i\in A'} q_{ki} \prod_{i\in B'} q_{i j(i)}=0$. Thus,
$R(Q)$ equals $ \prod_{i\in A} q_{ki} \prod_{i\in B} q_{i j(i)}$ times
the coefficient of the monomial. But $R(Q)\geq0$, so that this
coefficient is a positive integer
% $0$ or not less than $1$,
and (\ref{eq:der2}) follows.

%s4.2 #&#
\subsection{Proof of Theorem \texorpdfstring{\protect\ref{theo:order}}{2}}
\label{subsec:order}

Applying Theorem~\ref{theo:gene:finite}, we get that under the
assumptions of Theorem~\ref{theo:order}, there exists $K$ such that
\[
\Pb_{\theta_{0}} \bigl(\llVert f_{l,\theta}- f_{l,\theta_0} \rrVert
_{1} \leq2K w_{n}\vert Y_{1:n} \bigr)=1+
\mathrm{o}_{\Pb_{\theta_{0}}}(1).
\]
But if inequality \eqref{ineqsimple} holds, then as soon as
%
%e25 #&#
\begin{equation}
\label{l1order} \llVert f_{l,\theta}- f_{l,\theta_0} \rrVert
_{1} \lesssim w_{n}
\end{equation}
we get that, for any
$ j\in\{1,\ldots,k\}$, either
$\Pb_{\theta} (X_{1}=j ) \lesssim w_{n}$, or
\[
\exists i\in\{1,\ldots,k_{0}\},\qquad \Pb_{\theta}
(X_{1}=j )\bigl \| \gamma_{j}-\gamma_{i}^{0}
\bigr \|^{2} \lesssim w_{n}.
\]
%
% {\text{ and }} |\Pb_{\theta}\left(X_{1}=j\right)-\Pb_{\theta_0}
% \lesssim u_{n}.
Let us choose $\epsilon\leq\min_{i \neq j} \| \gamma_i^0 - \gamma_j^0
\|/4 $ in the definition of $B(i) $ in \eqref{ineqsimple}. We then
obtain that for large enough $n$, all $j_1, j_2 \in J(\theta)$, we have
$j_1 \sim j_2$ if and only if they belong to the same $B(i)$, $i =1,
\ldots, k_0$, so that $L(\theta) \leq k_0$. On the other hand,
$L(\theta) < k_0$ would mean that at least one $B(i)$ would be empty
which contradicts the fact that
\[
\bigl\llvert \Pb_{\theta} \bigl(X_{1}\in B(i) \bigr)-
\Pb_{\theta_0} (X_{1}=i )\bigr\rrvert \leq w_n.
\]
Thus, for large enough $n$, if \eqref{l1order} holds, then $L(\theta
)=k_0$, so that
\[
P^\Pi \bigl[ L(\theta)=k_{0} | Y^n \bigr] =1
+ \mathrm{o}_{\mathbb P_{\theta_0}}(1).
\]

To finish the proof, we now prove that \eqref{ineqsimple} holds under
the assumptions of Theorem~\ref{theo:order}. This will follow from
Proposition~\ref{lemcondi} below which is slightly more general.

An inequality that relates the $L_1$ distance of the $l$-marginals to
the parameters of the HMM
%, stated below as inequality (\ref{dul1aupara}),
is proved in Gassiat and van Handel \cite{HanGas} for translation mixture models, with the
strength of being uniform over the number (possibly infinite) of
populations in the mixture.
%
%we use \citet{HanGas}, where a theorem gives such transportation for
%any kind of mixture model.
However, for our purpose, we do not need such a general result, and it
is possible to obtain it
%prove that (\ref{dul1aupara}) holds
for more general situations than families of translated distributions,
under the structural assumption \textbf{A4}.
The inequality following Theorem~3.10 of
Gassiat and van Handel \cite{HanGas} says that
there exists a constant $c(\theta_{0})>0$ such that for any small
enough positive $\varepsilon$,
%
%e26 #&#
\begin{eqnarray}
\label{dul1aupara}
&&\frac{\llVert  f_{l,\theta}- f_{l,\theta_0} \rrVert _{1}}{c(\theta
_{0})}
\nonumber
\\
&&\quad\geq \sum_{1\leq j \leq k : \forall i, \|\gamma_{j}-\gamma_{i}^{0}\| >
\varepsilon}
\Pb_{\theta} (X_{1}=j )
\nonumber
\\
&&\qquad{}+ \sum_{1\leq i_{1},\ldots,i_{l} \leq k_{0}} \left[ \bigl\llvert \Pb_{\theta}
\bigl(X_{1:l}\in A(i_{1}, \ldots,i_{l}) \bigr)-
\Pb_{\theta
_0} (X_{1:l}=i_{1}\cdots i_{l} )
\bigr\rrvert\vphantom{\pmatrix{
\gamma_{j_{1}}
\cr
\cdots
\cr
\gamma_{j_{l}} }} \right.
\\
&&\phantom{\qquad{}+ \sum_{1\leq i_{1},\ldots,i_{l} \leq k_{0}} \Biggl[}{}+
 \left\llVert \sum_{( j_{1},\ldots,j_{l})\in A(i_{1},  \ldots,i_{l})} \Pb_{\theta}
(X_{1:l}=j_{1}\cdots j_{l} ) \left\{ \pmatrix{
\gamma_{j_{1}}
\cr
\cdots
\cr
\gamma_{j_{l}} } -\pmatrix{
\gamma_{i_{1}}^{0}
\cr
\cdots
\cr
\gamma_{i_{l}}^{0}
} \right\}\right\rrVert
\nonumber
\\
&& \left.\phantom{\qquad{}+ \sum_{1\leq i_{1},\ldots,i_{l} \leq k_{0}} \Biggl[}{}+
\frac{1}{2} \sum_{ ( j_{1},\ldots,j_{l})\in A(i_{1},  \ldots,i_{l})} \Pb_{\theta}
(X_{1:l}=j_{1}\cdots j_{l} )\left\llVert
\pmatrix{\gamma_{j_{1}}
\cr
\cdots
\cr
\gamma_{j_{l}} } - \pmatrix{
\gamma_{i_{1}}^{0}
\cr
\cdots
\cr
\gamma_{i_{l}}^{0}
}\right\rrVert ^{2} \right],
\nonumber
\end{eqnarray}
where $A(i_{1},  \ldots,i_{l})=\{( j_{1},\ldots,j_{l})\dvtx  \|\gamma
_{j_1}-\gamma_{i_1}^{0}\|\leq\varepsilon,\ldots, \|\gamma_{j_l}-\gamma
_{i_l}^{0}\|\leq\varepsilon\}$. The above lower bound essentially
corresponds to a partition of $\{ 1,\ldots, k\}^l$ into $k_0^l + 1$
groups, where the first $k_0^l$ groups correspond to the components
that are close to true distinct components in the multivariate mixture
and the last corresponds to components that are emptied. The first term
on the right-hand side controls the weights of the components that are
emptied (group $k_0^l + 1$), the second term controls the sum of the
weights of the components belonging to the $i$th group, for $i=1,\ldots
, k_0^l$ (components merging with the true $i$th component), the third
term controls the distance between the mean value over the group $i$
and the true value of the $i$th component in the true mixture while the
last term controls the distance between each parameter value in group
$i$ and the true value of the $i$th component.

Notice that \eqref{ineqsimple} is a consequence of \eqref{dul1aupara}.
We shall prove that \eqref{dul1aupara} holds under an assumption
slightly more general than \textbf{A4}.
For this, we need to introduce some notations.
For all $I = (i_1,\ldots, i_l) \in\{ 1,\ldots, k\}^l$,
% : = \mathcal I_l$ or eventually $I = (i_1,{\chr"C9},i_l) \in\{ 1,
define $\gamma_{I}=(\gamma_{i_{1}},\ldots,\gamma_{i_{l}})$,
% and note that
%$\gamma_{I}^0=(\gamma_{i_{1}},\ldots,\gamma_{i_{l}})$ necessarily
%corresponds to $I \in\mathcal I_{0,l}$. Write $ \mathcal I_{0,l}^c =
%I_{0,l}$ set $\mathcal J(I) = \{ 1,\cdots, t_{i_1}-t_{i_1-1}\}\times
%Define
$G_{\gamma_I} = \prod_{t = 1}^{l} g_{\gamma_{i_t}}(y_t)$, $D^1G_{\gamma
_I}$ the vector of first derivatives of $G_{\gamma_I}$ with respect to
each of
the
distinct
elements in $\gamma_I$, note that it has dimension $d \times| I |$,
where $|I| $ denotes the number of distinct indices in $I$, and
similarly define $D^2G_{\gamma_I}$ the symmetric matrix in $R^{d|I|
\times d |I|}$ made of second derivatives of $G_{\gamma_I}$ with
respect to the distinct elements (indices) in $\gamma_I$. For any
$\mathbf{t}=(t_{1},\ldots,t_{k_{0}})\in T$, define for all $i\in\{
1,\ldots,k_{0}\}$ the set
$J(i)=\{t_{i-1}+1,\ldots,t_{i}\}$, using $t_{0}=0$.

We then consider the following condition:
\begin{enumerate}[\textbf{A4bis}]
\item[\textbf{A4bis}]
For any $\mathbf{t}=(t_{1},\ldots,t_{k_{0}})\in T$,
for all collections $(\pi_I)_I$, $(\gamma_I)_I$, $I \!\notin\!\{1,\ldots
,t_{k_{0}}\}^l$ satisfying $\pi_I\geq0$,
$\gamma_{I}=(\gamma_{i_{1}},\ldots,\gamma_{i_{l}})$ such that $\gamma
_{i_j}=\gamma_{i}^{0}$ when $i_j\in J(i)$ for some $i\leq k_0$ and
$\gamma_{i_j}\in\Gamma\setminus\{ \gamma_i^0, i=1,\ldots,k_0\}$ when
$i_j\notin\{1,\ldots,t_{k_{0}}\}$, for all collections
$ ( a_I)_I$, $(c_I)_I$, $(b_I)_I$, $I\in\{1,\ldots,k_{0}\}^{l}$, $a_I\in\R
$, $c_I\geq0$ and $b_I \in\R^{d |I|}$, for all collection of vectors
$z_{I,J} \in\R^{d |I|}$ with $I\in\{1,\ldots,k_{0}\}^{l}$ and $J \in
J(i_{1})\times\cdots\times J(i_{l})$ satisfying $\|z_{I,J}\| = 1$,
and all sequences $(\alpha_{I,J})$, satisfying $\alpha_{I,J}\geq0$ and
$\sum_{J \in J(i_{1})\times\cdots\times J(i_{l})} \alpha_{I,J} =1$,
%%and for all collection of symmetric semi-definite matrices $M_I \in
%
%e27 #&#
\begin{eqnarray}
\label{iden:l} %
&& \sum_{I \notin\{1,\ldots,t_{k_{0}} \}^l}
\pi_I G_{\gamma_I} + \sum_{I\in\{1,\ldots,k_{0}\}^{l}}
\bigl( a_I G_{\gamma_I^0} + b_{I}^{T}
D^1G_{\gamma_I^0} \bigr)
\nonumber
\\
&&\qquad{} + \sum_{I\in \{1,\ldots,k_{0}\}^{l}} c_I \sum
_{J \in J(i_{1})\times
\cdots\times J(i_{l})} \alpha_{I,J}z_{I,J}^{T}
D^2 G_{\gamma_I^0} z_{I,J}= 0
\nonumber
\\[-8pt]
\\[-8pt]
%+ \mbox{tr}\left[ M_I D^2 G_{\gamma_I^0}\right]
&&\quad \Leftrightarrow
\nonumber
\\
&& \sum_{I \notin\{1,\ldots,t_{k_{0}} \}^l} \pi_I+ \sum
_{I\in \{
1,\ldots,k_{0}\}^{l}} \bigl( |a_I| + \|b_I\|
+ c_I \bigr) = 0.
\nonumber
%+
%
\end{eqnarray}
\end{enumerate}

We have the following proposition.
%
%pr1 #&#
\begin{prop}
\label{lemcondi}
Assume that the function $\gamma\mapsto g_{\gamma}(y)$ is twice
continuously differentiable in $\Gamma$ and that for all $y$, $g_{\gamma
}(y)$ vanishes as $\|\gamma\|$ tends to infinity. Then, if assumption
\textup{\textbf{A4bis}} is verified, \eqref{dul1aupara} holds. Moreover,
condition \textup{\textbf{A4bis}} is verified as soon as condition \textup{\textbf{A4}}
(corresponding to $l=1$) is verified.
\end{prop}

Let us now prove Proposition~\ref{lemcondi}.
To prove the first part of the proposition, we follow the ideas of the
beginning of the proof of Theorem~5.11 in Gassiat and van Handel \cite{HanGas}.
If \eqref{dul1aupara} does not hold, there exist a sequence of
$l$-marginals $(f_{l,\theta^n})_{n\geq1}$ with parameters $(\theta
^n)_{n\geq1}$ such that for some positive sequence $\varepsilon_{n}$
tending to $0$,
$\|f_{l,\theta^n}-f_{l,\theta_0}\|_{1} / N_{n} (\theta^n )$ tends to
$0$ as $n$ tends to infinity, with
\begin{eqnarray*}
N_{n} (\theta )&=& \sum_{1\leq j \leq l : \forall i, \|\gamma_{j}-\gamma_{i}^{0}\| >
\varepsilon_{n}}
\Pb_{\theta} (X_{1}=j )
\\
&&{}+ \sum_{1\leq i_{1},\ldots,i_{l} \leq k_{0}} \left[ \biggl|\sum
_{( j_{1},\ldots,j_{l})\in A_{n}(i_{1},  \ldots,i_{l})} \Pb_{\theta} (X_{1:l}=j_{1}
\cdots j_{l} )-\Pb_{\theta_0} (X_{1:l}=i_{1}
\cdots i_{l} )\biggr|\right.
\\
&&\phantom{{}+ \sum_{1\leq i_{1},\ldots,i_{l} \leq k_{0}} \Biggl[} {}+  \left
\llVert \sum_{( j_{1},\ldots,j_{l})\in A_{n}(i_{1},  \ldots,i_{l})} \Pb_{\theta}
(X_{1:l}=j_{1}\cdots j_{l} ) \left\{\pmatrix{
\gamma_{j_{1}}
\cr
\cdots
\cr
\gamma_{j_{l}} } -\pmatrix{
\gamma_{i_{1}}^{0}
\cr
\cdots
\cr
\gamma_{i_{l}}^{0}
} \right\}\right\rrVert
\\
&&\left.\phantom{{}+ \sum_{1\leq i_{1},\ldots,i_{l} \leq k_{0}} \Biggl[} {}+
\frac{1}{2} \sum_{ ( j_{1},\ldots,j_{l})\in A_{n}(i_{1},  \ldots,i_{l})} \Pb_{\theta}
(X_{1:l}=j_{1}\cdots j_{l} )\left\llVert
\pmatrix{\gamma_{j_{1}}
\cr
\cdots
\cr
\gamma_{j_{l}} } - \pmatrix{
\gamma_{i_{1}}^{0}
\cr
\cdots
\cr
\gamma_{i_{l}}^{0}
}\right\rrVert ^{2} \right]
\end{eqnarray*}
with $A_{n}(i_{1},  \ldots,i_{l})=\{( j_{1},\ldots,j_{l})\dvtx  \|\gamma
_{j_1}-\gamma_{i_1}^{0}\|\leq\varepsilon_{n},\ldots, \|\gamma
_{j_l}-\gamma_{i_l}^{0}\|\leq\varepsilon_{n} \}$.

Now,
$f_{l,\theta^n}=\sum_{I\in\{1,\ldots,k\}^l} \Pb_{\theta^n}
(X_{1},\ldots,X_{l}=I )G_{\gamma_{I}^{n}}$ where $\theta^n=(Q^n,
(\gamma_{1}^{n},\ldots,\gamma_{k}^{n}))$, $Q^{n}$ a transition matrix
on $\{1,\ldots,k\}$. It is possible to extract a subsequence along
which, for all $i=1,\ldots,k$, either $\gamma_{i}^{n}$ converges to
some limit $\gamma_{i}$ or $\|\gamma_{i}^{n}\|$ tends to infinity.
Choose now the indexation such that for $i=1,\ldots,t_{1}$, $\gamma
_{i}^{n}$ converges to $\gamma_{1}^{0}$, for $i=t_{1}+1,\ldots,t_{2}$,
$\gamma_{i}^{n}$ converges to $\gamma_{2}^{0}$, and so on,
for $i=t_{k_{0}-1}+1,\ldots,t_{k_{0}}$, $\gamma_{i}^{n}$ converges to
$\gamma_{k_{0}}^{0}$, and if $t_{k_{0}}<k$, for some $\tilde{k}\leq k$,
for $i=t_{k_{0}}+1,\ldots,\tilde{k}$, $\gamma_{i}^{n}$ converges to
some $\gamma_{i}\notin\{\gamma_{1}^{0},\ldots,\gamma_{k_{0}}^{0}\}$,
and for $i=\tilde{k}+1,\ldots,k$, $\|\gamma_{i}^{n}\|$ tends to
infinity. It is possible that $\tilde{k}=t_{k_{0}}$ in which case
no $\gamma_{i}^{n}$ converges to some $\gamma_{i}\notin\{\gamma
_{1}^{0},\ldots,\gamma_{k_{0}}^{0}\}$. Such a $\mathbf{t}=(t_{1},\ldots
,t_{k_{0}})\in T$ exists, because if $\|f_{l,\theta^n}-f_{l,\theta_0}\|
_{1} / N_{n} (\theta^n )$ tends to $0$ as $n$ tends to infinity, $\|
f_{l,\theta^n}-f_{l,\theta_0}\|_{1}$, and $N_{n} (\theta^n )$ tends to
$0$ as $n$ tends to infinity (if it was not the case, using the
regularity of $\theta\mapsto f_{l,\theta}$ we would have a
contradiction). Now along the subsequence we may write, for large
enough~$n$:
\begin{eqnarray*}
N_{n} \bigl(\theta^n \bigr)&=& \sum
_{I \notin\{1,\ldots,t_{k_{0}} \}^l} \Pb_{\theta} (X_{1:l}=I )
\\
&&{}+\sum_{I\in\{1,\ldots,k_{0}\}^{l}} \biggl[ \biggl\llvert \sum
_{J \in J(i_{1})\times\cdots\times J(i_{l})} \Pb_{\theta} (X_{1:l}=J )-
\Pb_{\theta_0} (X_{1:l}=I )\biggr\rrvert
\\
&&\phantom{+\sum_{I\in\{1,\ldots,k_{0}\}^{l}} \biggl[} {}  + \biggl\|\sum
_{J \in J(i_{1})\times\cdots\times J(i_{l})} \Pb_{\theta} (X_{1:l}=J )
\bigl(\gamma_{J} -\gamma_{I}^{0} \bigr)\biggr\|
\\
&&\phantom{+\sum_{I\in\{1,\ldots,k_{0}\}^{l}} \biggl[} {} +
\frac{1}{2} \sum_{J \in J(i_{1})\times\cdots\times J(i_{l})} \Pb_{\theta}
(X_{1:l}=J )\bigl\llVert \gamma_{J} -\gamma_{I}^{0}
\bigr\rrVert ^{2} \biggr].
\end{eqnarray*}
We shall use a Taylor expansion till order $2$.
To be perfectly rigorous in the following, we need to express
explicitly $I$ in terms of its distinct indices, $(\tilde{i}_1,\ldots,
\tilde{i}_{|I|})$, so that $G_{\gamma_I} = \prod_{t=1}^{|I|} \prod_{j:
i_j = \tilde{i}_t} g_{\gamma_{\,\tilde{i}_t}}(y_j)$, but to keep
notations concise we do not make such a distinction
and for instance
$ (\gamma_J^n-\gamma_I^0)^{T} D^1G_{\gamma_I^0}$ means\vspace*{-2pt}
\[
\sum_{t=1}^{|I|} \bigl(
\gamma_{\,\tilde{i}_t} -\gamma_{\,\tilde{i}_t}^0\bigr)^{T}
\frac
{ \partial G_{\gamma_I}}{ \partial\gamma_{\,\tilde{i}_t}},
\]
and similarly for the second derivatives. We have\vspace*{-2pt}
\begin{eqnarray*}
f_{l,\theta^n}-f_{l,\theta_0}& =& \sum_{I \notin\{1,\ldots,t_{k_{0}} \}^l}
\Pb_{\theta} (X_{1:l}=I ) G_{\gamma_{I}^n}
\\[-1pt]
&&{} +\sum_{I\in\{1,\ldots,k_{0}\}^{l}} \biggl\{ \biggl[\sum
_{J \in J(i_{1})\times\cdots\times J(i_{l})} \Pb_{\theta} (X_{1:l}=J )-
\Pb_{\theta_0} (X_{1:l}=I ) \biggr]G_{\gamma_{I}^0}
\\[-1pt]
&&\phantom{+\sum_{I\in\{1,\ldots,k_{0}\}^{l}} \biggl\{} {} + \sum
_{J \in J(i_{1})\times\cdots\times J(i_{l})} \Pb_{\theta} (X_{1:l}=J ) \bigl(
\gamma_{J} -\gamma_{I}^{0}
\bigr)^{T}D^1G_{\gamma_{I}^0}
\\[-1pt]
&&\phantom{+\sum_{I\in\{1,\ldots,k_{0}\}^{l}} \biggl\{} {}   +
\frac{1}{2} \sum_{J \in J(i_{1})\times\cdots\times J(i_{l})} \Pb_{\theta}
(X_{1:l}=J ) \bigl(\gamma_{J} -\gamma_{I}^{0}
\bigr)^{T}D^2G_{\gamma_{I}^*} \bigl(\gamma_{J}
-\gamma_{I}^{0} \bigr) \biggr\}
\end{eqnarray*}
with $ \gamma_I^*\in(\gamma_{I}^n, \gamma_I^0)$.
Thus, using the fact that for all $y$, $g_{\gamma}(y)$ vanishes as $\|
\gamma\|$ tends to infinity, $f_{l,\theta^n}-f_{l,\theta_0}/ N_{n}
(\theta^n )$ converges pointwise along a subsequence to a function $h$
of form\vspace*{-1pt}
\begin{eqnarray*}
h&=& \sum_{I \notin\{1,\ldots,t_{k_{0}} \}^l} \pi_I
G_{\gamma_I} + \sum_{I\in\{1,\ldots,k_{0}\}^{l}} \bigl(
a_I G_{\gamma_I^0} + b_{I}^{T}
D^1G_{\gamma_I^0} \bigr)
\\[-1pt]
&&{}+ \sum_{I\in \{1,\ldots,k_{0}\}^{l}} c_I \sum
_{J \in J(i_{1})\times
\cdots\times J(i_{l})} \alpha_{I,J}z_{I,J}^{T}
D^2 G_{\gamma_I^0} z_{I,J}
\end{eqnarray*}
as in condition $L(l)$, with $\sum_{I \notin\{1,\ldots,t_{k_{0}} \}^l}
\pi_I+ \sum_{I\in \{1,\ldots,k_{0}\}^{l}}  ( |a_I| + \|b_I\| + c_I
)=1$.
But as $\|f_{l,\theta^n}-f_{l,\theta_0}\|_{1} / N_{n} (\theta^n )$
tends to $0$ as $n$ tends to infinity, we have $\|h\|_{1} =0$ by
Fatou's lemma, and thus $h=0$, contradicting the assumption.

Let us now prove that \textbf{A4} implies \textbf{A4bis}. Let\vspace*{-1pt}
\begin{eqnarray*}
&& \sum_{I \notin\{1,\ldots,t_{k_{0}} \}^l} \pi_I G_{\gamma_I}
+ \sum_{I\in\{1,\ldots,k_{0}\}^{l}} \bigl( a_I
G_{\gamma_I^0} + b_{I}^{T} D^1G_{\gamma_I^0}
\bigr)
\\[-1pt]
&&\qquad{} + \sum_{I\in \{1,\ldots,k_{0}\}^{l}} c_I \sum
_{J \in J(i_{1})\times
\cdots\times J(i_{l})} \alpha_{I,J}z_{I,J}^{T}
D^2 G_{\gamma_I^0} z_{I,J}= 0
\\[-1pt]
%+ \mbox{tr}\left[ M_I D^2 G_{\gamma_I^0}\right]
&&\quad \Leftrightarrow
\\[-1pt]
&& \sum_{I \notin\{1,\ldots,t_{k_{0}} \}^l} \pi_I+ \sum
_{I\in \{
1,\ldots,k_{0}\}^{l}} \bigl( |a_I| + \|b_I\|
+ c_I \bigr) = 0 %+
\end{eqnarray*}
with $\pi_I$, $a_I$, $b_I$, $\alpha_{I,J} $ and $z_{I,J}$ be as in assumption
\textbf{A4bis}.
We group the terms depending only on $y_1$ and we can rewrite the
equation as
%
%e28 #&#
\begin{eqnarray}
\label{hy1} %
&& \sum_{i=t_{k_0}+1}^k
\pi_i'(y_2,\ldots,y_{l})g_{\gamma_i}(y_1)
+ \sum_{i=1}^{k_0} \bigl(
a_i'(y_2,\ldots, y_{l})
g_{\gamma_i^0}(y_1) + b_{i}^{\prime T}(y_2,
\ldots, y_{l}) D^1g_{\gamma_i^0}(y_1)
\bigr)
\nonumber
\\[-8pt]
\\[-8pt]
&&\quad{} +\sum_{i=1}^{k_{0}}\sum
_{j=1}^{t_i-t_{i-1}} \sum_{i_2,\ldots,i_l =1}^{k_0}
c_I' \sum_{(j_2,\ldots,j_l) \in
J(i_2)\times\cdots\times J(i_l)}
\alpha_{I,J} z_{I,J}(i)^{T} D^2g_{\gamma_i^0}(y_1)
z_{I,J}(i) =0,
\nonumber
\end{eqnarray}
where we have written
\[
z_{I,J} = \bigl(z_{I,J}(i_1),
\ldots,z_{I,J}(i_l)\bigr), \qquad\mbox{with } I =
(i,i_2,\ldots,i_l), J = (j_1,
\ldots,j_l), z_{I,J}(i) \in\R^d
\]
and
\[
c_I' = c_I \prod
_{t=2}^l g_{\gamma_{i_t}^0}(y_t).
\]
Note that if for $i=1,\ldots,k_0$ and $j = 1, \ldots, t_i - t_{i-1}$,
there exists $w_{i,j} \in\R^d$ such that
\[
\sum_{i_2,\ldots,i_l =1}^{k_0} c_I'
\sum_{(j_2,\ldots,j_l) \in
J(i_2)\times\cdots\times J(i_l)} \alpha_{I,J}
z_{I,J}(i)^{T} D^2g_{\gamma_i^0}(y_1)
z_{I,J}(i) = w_{i,j}^{T} D^2g_{\gamma_i^0}(y_1)
w_{i,j},
\]
where possibly $w_{i,j}= 0$. Let $\alpha_{i,j} = \|w_{i,j}\|^2/( \sum_{j=1}^{t_i-t_{i-1}} \|w_{i,j}\|^2)$ if there exists $j$ such that $\|
w_{i,j}\|^2 >0$
and $c_i' = \sum_{i_2,\ldots, i_l} c_I' \sum_{j=1}^{t_i-t_{i-1}} \|
w_{i,j}\|^2$, then
\begin{eqnarray*}
&&\sum_{j=1}^{t_i- t_{i-1}}\sum
_{i_2,\ldots,i_l =1}^{k_0} c_I'\sum
_{(j_2,\ldots,j_l) \in J(i_2)\times\cdots\times J(i_l)} \alpha_{I,J} z_{I,J}(i)^{T}
D^2g_{\gamma_i^0}(y_1) z_{I,J}(i)
\\
&&\quad= c_i' \sum_{j=1}^{t_i- t_{i-1}}
\alpha_{i,j} w_{i,j}^{T} D^2g_{\gamma_i^0}(y_1)
w_{i,j}
\end{eqnarray*}
and \eqref{iden:1} implies that
\[
a_i' = c_i' = 0,\qquad
b_i'= 0,\qquad i=1,\ldots,k_0, \qquad
\pi_i' = 0,\qquad i = t_{k_0}+1, \ldots, k.
\]
Simple calculations imply that
\[
\pi_i' = \sum_{i_2,\ldots, i_l= 1}^{k}
\pi_I \prod_{t=2}^{l}
g_{\gamma
_{i_t}^0}(y_t) = 0 \quad\Leftrightarrow\quad
\forall(i_2,\ldots,i_l) \in\{ 1,\ldots, k
\}^{l-2} \pi_{i,i_2,\ldots, i_l}=0
\]
and similarly if $i$ is such that there exists $j=1,\ldots,t_i -
t_{i-1}$, $I = (i,i_2,\ldots,i_l)$ and $J = (j,j_2,\ldots, j_l) \in
J(i)\times\cdots\times J(i_l)$ such that $c_I>0$, $\alpha_J >0$ and
$\| z_{I,J}(i)\| >0$, then $c_{i,i_2,\ldots, i_l}= 0$ for all
$i_2,\ldots, i_l$. Else, by considering $y_t$ for some other $t$, we
obtain that
\eqref{hy1} implies that
\[
\pi_I = 0 \qquad\forall I \notin\{1,\ldots, t_{k_0}
\}^l,\qquad c_I= 0 \qquad\forall I \in\{1,\ldots,
t_{k_0}\}^l.
\]
This leads to
\[
b_i' = \sum_{i_2,\ldots, i_l=1}^{k_0}
b_I \prod_{t\geq2}g_{\gamma
_{i_t}^0}(y_t)
= 0 \qquad\forall i =1, \ldots, k_0.
\]
A simple recursive argument implies that $b_I = 0$ for all $I \in\{
1,\ldots, t_{k_0}\}^l$ which in turns implies that $a_I= 0$
for all $I \in\{1,\ldots, t_{k_0}\}^l$ and condition \textbf{A4bis}
is verified.

%s4.3 #&#
\subsection{Proof of Theorem \texorpdfstring{\protect\ref{theo:uncontredeux}}{3}}
\label{subsec:uncontredeux}

First, we obtain
the following lemma.
%
%le2 #&#
\begin{lemma}
\label{lemrootn}
Under the assumptions of Theorem~\ref{theo:uncontredeux}, for any
sequence $M_{n}$ tending to infinity,
\[
\Pb^{\Pi} \biggl( (p+q ) \wedge \bigl(2-(p+q) \bigr)\llVert
f_{2,\theta}- f_{2,\theta_0} \rrVert _{1} \leq
\frac{M_{n}}{\sqrt {n}} \biggr)=1+\mathrm{o}_{\Pb_{\theta_0}}(1).
\]
\end{lemma}

We prove Lemma~\ref{lemrootn} by applying Theorem~\ref{theo:gene},
using some of the computations of the proof of
Theorem~\ref{theo:gene:finite} but verifying assumption \textbf{C3bis} instead of
\textbf{C3}.
Set $S_{n}=U_{n}\times\mathcal{X}$ with
\begin{eqnarray*}
U_{n}&=& \biggl\{\theta=(p,q,\gamma_{1},\gamma_{2})
\dvtx \bigl \|\gamma_{1}-\gamma ^{0}\bigr \|^{2}\leq
\frac{1}{\sqrt{n}}, \bigl \|\gamma_{2}-\gamma^{0}\bigr \|^{2}
\leq \frac{1}{\sqrt{n}},
\\
&&\phantom{\biggl\{} {}\bigl \|q\bigl(\gamma_{1}-
\gamma^{0}\bigr)+p\bigl(\gamma_{2}-\gamma^{0}
\bigr)\bigr \|\leq\frac{1}{\sqrt {n}}, \biggl|q-\frac{1}{2}\biggr|\leq\epsilon, \biggl|p-
\frac{1}{2}\biggr|\leq\epsilon  \biggr\}
\end{eqnarray*}
for small but fixed $\epsilon$. We shall prove later the following lemma.

%le3 #&#
\begin{lemma} \label{k2k01:Dn}
Let $M_{n}$ tend to infinity. Then
\[
\sup_{(\theta, x)\in S_n} \Pb_{\theta_0} \bigl[ \ell_n(
\theta,x) - \ell _n(\theta_0, x_0) <
-M_n \bigr] = \mathrm{o}(1)
\]
and
%
%e29 #&#
\begin{equation}
\label{Sn:D} \Pi(S_n) \gtrsim n^{-3d/4}.
\end{equation}
\end{lemma}

% From Lemma~\ref{k2k01:Dn} we obtain that \eqref{denom} holds with
%such a definition of $S_n$ and that
%(instead of $n^{-d}$ obtained in the proof of Theorem
Now we prove that assumption \textbf{C3bis} holds with $\epsilon_n =
M_n /\sqrt{n}$, which will finish the proof of Lemma~\ref{lemrootn}.
%
%if \eqref{denom} and \eqref{Sn:D} hold, which will be proved later.
%to eliminate the term $\sqrt{\log n}$ in the concentration rate of $
By Proposition~\ref{lemcondi},
we obtain that there exists $c(\theta_{0})>0$ and $\eta>0$ such that:
\begin{itemize}
\item
If $\|\gamma_{1}-\gamma^{0}\|\leq\eta$ and $\|\gamma_{2}-\gamma^{0}\|
\leq\eta$,
\begin{eqnarray*}
&&\llVert f_{2,\theta}- f_{2,\theta_0} \rrVert _{1}
\\
&&\quad\geq c(\theta _{0})\frac{1}{p+q} \bigl[\bigl\llVert q
\bigl(\gamma_{1}-\gamma^{0}\bigr)+p\bigl(
\gamma_{2}-\gamma^{0}\bigr)\bigr\rrVert +q\bigl\llVert
\gamma_{1}-\gamma^{0}\bigr\rrVert ^{2} + p \bigl
\llVert \gamma_{2}-\gamma^{0}\bigr\rrVert ^{2}
\bigr].
\end{eqnarray*}
\item
If $\|\gamma_{1}-\gamma^{0}\|\leq\eta$ and $\|\gamma_{1}-\gamma^{0}\|
+\|\gamma_{2}-\gamma^{0}\| > 2\eta$,
\[
\llVert f_{2,\theta}- f_{2,\theta_0} \rrVert _{1} \geq c(
\theta_{0}) \biggl[\frac{p}{p+q} +\frac{q}{p+q} \bigl\llVert
\gamma_{1}-\gamma^{0}\bigr\rrVert \biggr].
\]
\item
If $\|\gamma_{2}-\gamma^{0}\|\leq\eta$ and $\|\gamma_{1}-\gamma^{0}\|
+\|\gamma_{2}-\gamma^{0}\| >2 \eta$,
\[
\llVert f_{2,\theta}- f_{2,\theta_0} \rrVert _{1} \geq c(
\theta_{0}) \biggl[\frac{q}{p+q} +\frac{p}{p+q} \bigl\llVert
\gamma_{2}-\gamma^{0}\bigr\rrVert \biggr].
\]
\item
If $\|\gamma_{1}-\gamma^{0}\|> \eta$ and $\|\gamma_{2}-\gamma^{0}\|
>\eta$,
\[
\llVert f_{2,\theta}- f_{2,\theta_0} \rrVert _{1} \geq c(
\theta_{0}).
\]
\end{itemize}
Similar upper bounds hold also by Taylor expansion. Thus, for any $m$,
$A_{n,m}(\epsilon_{n})$ is a subset of the set of $\theta$'s such that
\begin{eqnarray*}
&&\min \biggl\{ \frac{(p+q)\wedge (2-(p+q))}{p+q} \bigl[\bigl\llVert q\bigl(\gamma_{1}-
\gamma^{0}\bigr)+p\bigl(\gamma_{2}-\gamma^{0}
\bigr)\bigr\rrVert +q\bigl\llVert \gamma_{1}-\gamma^{0}\bigr
\rrVert ^{2} + p \bigl\llVert \gamma_{2}-
\gamma^{0}\bigr\rrVert ^{2} \bigr];
\\
&&\phantom{\min\biggl\{} \frac{(p+q)\wedge (2-(p+q))}{p+q} \bigl[p +q \bigl\llVert
\gamma_{1}-\gamma^{0}\bigr\rrVert \bigr];
\\
&&\phantom{\min\biggl\{} \frac{(p+q)\wedge (2-(p+q))}{p+q} \bigl[q +p \bigl\llVert
\gamma_{2}-\gamma^{0}\bigr\rrVert \bigr]; (p+q)\wedge
\bigl(2-(p+q)\bigr) \biggr\} \lesssim(m+1)\epsilon_{n}.
\end{eqnarray*}
This leads to
\[
\Pi_{2} \bigl(A_{n,m}(\epsilon_{n}) \bigr)
\lesssim\bigl[(m+1)\epsilon_{n}\bigr]^{2\alpha}+\bigl[(m+1)
\epsilon_{n}\bigr]^{2\beta
}+\bigl[(m+1)\epsilon_{n}
\bigr]^{\alpha+d}
\]
so that if $\alpha, \beta>3d/4$ and \eqref{Sn:D} holds, there exists
$\delta>0$ such that
\[
\frac{\Pi_{2}  (A_{n,m}(\epsilon_{n}) ) \mathrm{e}^{- \trup{( n
m^2\epsilon_n^2)}{( 32 l)}}}{ \Pi(S_n) } \lesssim n^{-\delta}\bigl[(M_nm)^{2\alpha}+(M_nm)^{2\beta}
+ (M_nm)^{\alpha
+d}\bigr]\mathrm{e}^{- \trup{( M_n^2 m^2 )}{( 32 l)} }.
\]
Also for all $\epsilon>0$ small enough $A_{n,m}(\epsilon)$ contains the
set of $\theta$'s such that
\begin{eqnarray*}
&&\max \biggl\{ \frac{(p+q)\wedge (2-(p+q))}{p+q}
\\
&&\phantom{\max\biggl\{}{}\times\bigl[\bigl\llVert q\bigl(\gamma_{1}-
\gamma^{0}\bigr)+p\bigl(\gamma_{2}-\gamma^{0}
\bigr)\bigr\rrVert +q\bigl\llVert \gamma_{1}-\gamma^{0}\bigr
\rrVert ^{2} + p \bigl\llVert \gamma_{2}-
\gamma^{0}\bigr\rrVert ^{2} \bigr];
\\
&&\phantom{\max\biggl\{} \frac{(p+q)\wedge (2-(p+q))}{p+q} \bigl[p +q \bigl\llVert
\gamma_{1}-\gamma^{0}\bigr\rrVert \bigr];
\\
&&\phantom{\max\biggl\{}\frac{(p+q)\wedge (2-(p+q))}{p+q} \bigl[q +p \bigl\llVert
\gamma_{2}-\gamma^{0}\bigr\rrVert \bigr]; (p+q)\wedge
\bigl(2-(p+q)\bigr)  \biggr\} \lesssim(m+1)\epsilon
\end{eqnarray*}
therefore
\[
N\biggl(\frac{m\epsilon_{n}}{12}, A_{n,m} (\epsilon_{n} ),
d_{l}(\cdot,\cdot) \biggr) \lesssim m^{2+2d} \lesssim
\mathrm{e}^{ \trup{( n\epsilon
_n^2 m^2(\rho_{\theta_0}-1)^2)}{( 16 l (2+ \rho_{\theta_0}-1)^2)}},
\]
so that assumption \textbf{C3bis} is verified.

We now prove Theorem~\ref{theo:uncontredeux}. Notice first that, by
setting
\[
D_{n}= \int_{\Theta\times{\mathcal X}} \mathrm{e}^{\ell_n(\theta,x) -
\ell_n(\theta_0,x_{0})}\Pi_{2}(\mathrm{d}\theta)\pi_{\mathcal X} (\mathrm{d}x
),
\]
 as in the proof of Theorem~\ref{theo:gene} we get that for any
sequence $C_{n}$ tending to infinity,
%
%e30 #&#
\begin{equation}
\label{denom3/4} \Pb_{\theta_{0}} \bigl(D_{n} \leq C_n
n^{-D/2} \bigr)=\mathrm{o} (1 )
\end{equation}
with $D=d+d/2$.

Let now $\epsilon_{n}$ be any sequence going to $0$
and let
$A_{n}= \{\frac{p}{p+q} \leq\epsilon_{n} \mbox{ or } \frac
{q}{p+q} \leq\epsilon_{n}
\}$.
For some sequence $M_{n}$ going to infinity and $\delta_{n}=M_{n}/\sqrt {n}$, let
$B_{n}=\{ (p+q ) \wedge (2-(p+q) )\llVert  f_{2,\theta
}- f_{2,\theta_0} \rrVert _{1} \leq\delta_{n}\}$.
We then control with $D=d+d/2$, using Lemma~\ref{lemrootn}
\begin{eqnarray*}
E_{\theta_{0}} \bigl[\Pb^{\Pi} (A_{n}\vert
Y_{1:n} ) \bigr]&=&E_{\theta_{0}} \bigl[\Pb^{\Pi}
(A_{n}\cap B_{n}\vert Y_{1:n} ) \bigr]+\mathrm{o}
(1 )
\\
&=&E_{\theta_{0}} \biggl[\frac{ \int_{A_n \cap B_{n} \times{\mathcal X}
} \mathrm{e}^{\ell_n(\theta,x) -\ell_n(\theta_0,x_{0})}\Pi_{2}(\mathrm{d}\theta) \pi
_{\mathcal X} (\mathrm{d}x ) }{ \int_{\Theta\times{\mathcal X}} \mathrm{e}^{\ell
_n(\theta,x) - \ell_n(\theta_0,x_{0})}\Pi_{2}(\mathrm{d}\theta) \pi_{\mathcal
X} (\mathrm{d}x ) } \biggr]+\mathrm{o} (1 )
\\
&:=&E_{\theta_{0}} \biggl[\frac{N_n}{D_n} \biggr]+\mathrm{o} (1 )
\\
&\leq& \Pb_{\theta_{0}} \bigl(D_{n} \leq C_n
n^{-D/2} \bigr)+ \frac
{n^{D/2}}{C_n}\Pi_{2}(A_{n}
\cap B_{n})+\mathrm{o} (1 ).
\end{eqnarray*}
Thus using \eqref{denom3/4}, the first part of Theorem~\ref{theo:uncontredeux} is proved by showing that
%
%e31 #&#
\begin{equation}
\label{eq:priorvider} \Pi_{2}(A_{n} \cap B_{n})
\lesssim\delta_{n}^{2\alpha}+\delta _{n}^{\alpha+d}+
\delta_{n}^{d+d/2}\epsilon_{n}^{\alpha-d/2}.
\end{equation}
Then, the second part of Theorem~\ref{theo:uncontredeux} follows from
its first part and Lemma~\ref{lemrootn}.

We now prove that (\ref{eq:priorvider}) holds.
Define
\begin{eqnarray*}
B_{n}^{1}&=& \biggl\{ \frac{(p+q)\wedge (2-(p+q))}{p+q}
\\
&&\phantom{\biggl\{} {}\times \bigl[\bigl\llVert q\bigl(\gamma_{1}-
\gamma^{0}\bigr)+p\bigl(\gamma_{2}-\gamma^{0}
\bigr)\bigr\rrVert +q\bigl\llVert \gamma_{1}-\gamma^{0}\bigr
\rrVert ^{2} + p \bigl\llVert \gamma_{2}-
\gamma^{0}\bigr\rrVert ^{2} \bigr]\leq\delta_{n}
  \biggr\},
\\
B_{n}^{2}&=& \biggl\{ \frac{(p+q)\wedge (2-(p+q))}{p+q} \bigl[p +q \bigl
\llVert \gamma_{1}-\gamma^{0}\bigr\rrVert \bigr] \leq
\delta_{n} \biggr\},
\\
B_{n}^{3}&=& \biggl\{ \frac{(p+q)\wedge (2-(p+q))}{p+q} \bigl[q +p \bigl
\llVert \gamma_{2}-\gamma^{0}\bigr\rrVert \bigr] \leq
\delta_{n} \biggr\}
\end{eqnarray*}
and
\[
B_{n}^{4}= \bigl\{ (p+q)\wedge \bigl(2-(p+q)\bigr)\leq
\delta_{n} \bigr\}.
\]
Then
\[
\Pi_{2}(A_{n} \cap B_{n}) \leq
\Pi_{2}\bigl(A_{n} \cap B_{n}^{1}
\bigr)+\Pi _{2}\bigl(A_{n} \cap B_{n}^{2}
\bigr)+\Pi_{2}\bigl(A_{n} \cap B_{n}^{3}
\bigr)+\Pi_{2}\bigl(A_{n} \cap B_{n}^{4}
\bigr).
\]
Notice that on $A_{n}$, if $p+q \geq1$, then $p\leq\epsilon_{n}$ and
$q\geq1-\epsilon_{n}$, or $q\leq\epsilon_{n}$ and $p\geq1-\epsilon
_{n}$, so that also
$2-(p+q)\geq1-\epsilon_{n}$. %
\begin{itemize}
\item
On $A_{n} \cap B_{n}^{1}$, $\llVert  q(\gamma_{1}-\gamma^{0})+p(\gamma
_{2}-\gamma^{0})\rrVert \lesssim\delta_{n}$,
$q\llVert  \gamma_{1}-\gamma^{0}\rrVert ^{2}\lesssim\delta_{n}$, $ p
\llVert  \gamma_{2}-\gamma^{0}\rrVert ^{2}\lesssim\delta_{n}$, and
$p\lesssim\epsilon_{n}$ or $q\lesssim\epsilon_{n}$. This gives $\Pi
_{2}(A_{n} \cap B_{n}^{1})\lesssim\delta_{n}^{d+d/2}\epsilon
_{n}^{\alpha- d/2}$.
\item
On $A_{n} \cap B_{n}^{2}$,
$p\lesssim\delta_{n}$ and $q\llVert  \gamma_{1}-\gamma^{0}\rrVert
\lesssim\delta_{n}$ in case $p+q\leq1$, and
$p\lesssim\delta_{n}$, $1-q\lesssim\delta_{n}$ and $q\llVert  \gamma
_{1}-\gamma^{0}\rrVert \lesssim\delta_{n}$ in case $p+q\geq1$,
leading to $\Pi_{2}(A_{n} \cap B_{n}^{2})\lesssim\delta_{n}^{\alpha+
d} +\delta_{n}^{\alpha+ \beta+ d} $.
\item
For symmetry reasons, $\Pi_{2}(A_{n} \cap B_{n}^{3})=\Pi_{2}(A_{n} \cap
B_{n}^{2})$.
\item
On $A_{n} \cap B_{n}^{4}$, $p\lesssim\delta_{n}$ and $q\lesssim\delta
_{n}$, so that $\Pi_{2}(A_{n} \cap B_{n}^{4})\lesssim\delta
_{n}^{2\alpha}$.
\end{itemize}
Keeping only the leading terms, we see that (\ref{eq:priorvider}) holds
and this terminates the proof Theorem~\ref{theo:uncontredeux}.

We now prove Lemma~\ref{k2k01:Dn}.
%, which will finish the proof of Theorem~\ref{theo:uncontredeux}.
%Notice that $U_{n}\subset B_{n}$, and
We easily get $\Pi_{2}(U_{n}) \gtrsim n^{-3d/4}$, and
\[
D_{n}\geq\int_{U_{n}\times{\mathcal X}} \mathrm{e}^{\ell_n(\theta,x) - \ell
_n(\theta_0,x)}
\Pi_{2}(\mathrm{d}\theta) \pi_{\mathcal X} (\mathrm{d}x ).
\]
Let us now study $\ell_n(\theta,x) - \ell_n(\theta_0,x)$. First,
following the proof of Lemma~2 of Douc \textit{et al.}
\cite{douc:moulines:ryden:04} we
find that, for any $\theta\in U_{n}$, for any $x$,
\[
\bigl\llvert \ell_{n} (\theta )-\ell_{n} (\theta,x )\bigr
\rrvert \leq \biggl(\frac{1+2\epsilon}{1-2\epsilon} \biggr)^{2},
\]
where $ \ell_{n} (\theta ) = \sum_{x=1}^k \mu_\theta(x) \ell
_n(\theta, x)$.
Thus, for any $\theta\in U_{n}$ and any $x$, and since $\ell_n(\theta
_0,x)$ does not depend on $x$,
%
%e32 #&#
\begin{equation}
\label{mindeno1} \ell_n(\theta,x) - \ell_n(
\theta_0,x)\geq\ell_{n}(\theta)-\ell_n(\theta
_0)- \biggl(\frac{1+2\epsilon}{1-2\epsilon} \biggr)^{2}.
\end{equation}
Let us now study $\ell_{n}(\theta)-\ell_n(\theta_0)$.
%Note $g_{\gamma}(\cdot)=g(\cdot- \gamma)$.
\begin{eqnarray*}
&&\ell_n(\theta) - \ell_n(\theta_0)
\\
&&\quad=\sum_{k=1}^{n}\log \biggl[
\Pb_{\theta
} (X_k=1 \vert Y_{1:k-1} )
\frac{g_{\gamma_{1}}}{g_{\gamma
^{0}}} (Y_{k} ) + \Pb_{\theta} (X_k=2
\vert Y_{1:k-1} ) \frac{g_{\gamma_{2}}}{g_{\gamma^{0}}} (Y_{k} ) \biggr]
\end{eqnarray*}
and we set for $k=1$
\begin{eqnarray*}
\Pb_{\theta} (X_k=1 \vert Y_{1:k-1} )=
\Pb_{\theta} (X_1=1 )&=&\frac{q}{p+q},
\\
\Pb_{\theta} (X_k=2 \vert Y_{1:k-1} )=
\Pb_{\theta} (X_1=2 )&=&\frac{p}{p+q}.
\end{eqnarray*}
Denote $p_{k}(\theta)$ the random variable $\Pb_{\theta} (X_k=1
\vert Y_{1:k-1}  )$, which is a function of $Y_{1:k-1}$ and thus
independent of $Y_{k}$. We have the recursion
%
%e33 #&#
\begin{equation}
\label{recursion} p_{k+1} (\theta )=\frac{(1-p)p_{k}(\theta)g_{\gamma
_{1}}(Y_{k})+q(1-p_{k}(\theta))g_{\gamma_{2}}(Y_{k})}{p_{k}(\theta
)g_{\gamma_{1}}(Y_{k})+(1-p_{k}(\theta))g_{\gamma_{2}}(Y_{k})}.
\end{equation}
Note that, for any $p$, $q$ in $]0,1[$, for any $k\geq1$,
\[
p_{k}\bigl(p,q,\gamma^{0},\gamma^{0}\bigr)=
\frac{q}{p+q}.
\]
We shall denote by $D^{i}_{(\gamma_{1})^{j},(\gamma_{2})^{i-j}}$ the
$i$th partial derivative operator $j$ times with respect to $\gamma
_{1}$ and $i-j$ times with respect to $\gamma_{2}$ ($0\leq j\leq i$,
the order in which derivatives are taken does not matter). Fix $\theta
=(p,q,\gamma_{1},\gamma_{2})\in U_{n}$. When derivatives are taken at point
$(p,q,\gamma^{0},\gamma^{0})$, they are written with $0$ as
superscript.

Using Taylor expansion till order $4$, there exists $t\in[0,1]$ such
that denoting $\theta_{t}=t\theta+ (1-t) (p,q,\gamma^{0},\gamma^{0})$:
%
%e34 #&#
\begin{equation}
\label{Taylor1} \ell_n(\theta) - \ell_n(
\theta_0)=\bigl(\gamma_{1}-\gamma^{0}
\bigr)D^{1}_{\gamma
_{1}}\ell_{n}^{0}+\bigl(
\gamma_{2}-\gamma^{0}\bigr)D^{1}_{\gamma_{2}}
\ell_{n}^{0} +S_{n}(\theta)+T_{n}(
\theta) + R_{n}(\theta,t),
\end{equation}
where $S_{n}(\theta)$ denotes the term of order $2$, $T_{n}(\theta)$
denotes the term of order $3$, and $R_{n}(\theta,t)$ the remainder,
that is
\begin{eqnarray*}
S_{n}(\theta)&=&\bigl(\gamma_{1}-\gamma^{0}
\bigr)^{2}D^{2}_{(\gamma_{1})^{2}}\ell _{n}^{0}+2
\bigl(\gamma_{1}-\gamma^{0}\bigr) \bigl(
\gamma_{2}-\gamma^{0}\bigr)D^{2}_{\gamma
_{1},\gamma_{2}}
\ell_{n}^{0} +\bigl(\gamma_{2}-
\gamma^{0}\bigr)^{2}D^{2}_{(\gamma_{2})^{2}}
\ell_{n}^{0},
\\
T_{n}(\theta)&=&\bigl(\gamma_{1}-\gamma^{0}
\bigr)^{3}D^{3}_{(\gamma_{1})^{3}}\ell _{n}^{0}+3
\bigl(\gamma_{1}-\gamma^{0}\bigr)^{2}\bigl(
\gamma_{2}-\gamma ^{0}\bigr)D^{3}_{(\gamma_{1})^{2},\gamma_{2}}
\ell_{n}^{0}
\\
&&{}+3\bigl(\gamma_{1}-\gamma^{0}\bigr) \bigl(
\gamma_{2}-\gamma^{0}\bigr)^{2}D^{3}_{\gamma
_{1},(\gamma_{2})^{2}}
\ell_{n}^{0} +\bigl(\gamma_{2}-
\gamma^{0}\bigr)^{3}D^{3}_{(\gamma_{2})^{3}}
\ell_{n}^{0}
\end{eqnarray*}
and
\[
R_{n}(\theta,t)=\sum_{k=0}^{4}
\pmatrix{k
\cr
4 }\bigl(\gamma_{1}-\gamma^{0}
\bigr)^{k}\bigl(\gamma_{2}-\gamma ^{0}
\bigr)^{4-k}D^{4}_{(\gamma_{1})^{k},(\gamma_{2})^{4-k}}\ell_n(
\theta_{t}).
\]
Easy but tedious computations lead to the following results.
\begin{eqnarray*}
 \bigl(\gamma_{1}-\gamma^{0}\bigr)D^{1}_{\gamma_{1}}
\ell_{n}^{0}+\bigl(\gamma_{2}-\gamma
^{0}\bigr)D^{1}_{\gamma_{2}}\ell_{n}^{0}
&=& \Biggl[\sum_{k=1}^{n}
\frac{D^{1}_{\gamma}g_{\gamma_{0}}}{g_{\gamma
_{0}}} (Y_{k} ) \Biggr] \biggl[\frac{q(\gamma_{1}-\gamma
^{0})+p(\gamma_{2}-\gamma^{0})}{p+q} \biggr]
\\
&=& \Biggl[\frac{1}{\sqrt{n}}\sum_{k=1}^{n}
\frac{D^{1}_{\gamma}g_{\gamma
_{0}}}{g_{\gamma_{0}}} (Y_{k} ) \Biggr] \biggl[\sqrt{n}\frac
{q(\gamma_{1}-\gamma^{0})+p(\gamma_{2}-\gamma^{0})}{p+q}
\biggr]
\end{eqnarray*}
so that
%
%e35 #&#
\begin{equation}
\label{ordre1} \sup_{\theta\in U_{n}}\bigl\llvert \bigl(
\gamma_{1}-\gamma^{0}\bigr)D^{1}_{\gamma
_{1}}
\ell_{n}^{0}+\bigl(\gamma_{2}-
\gamma^{0}\bigr)D^{1}_{\gamma_{2}}\ell_{n}^{0}
\bigr\rrvert =\mathrm{O}_{\Pb_{\theta_{0}}} (1 ).
\end{equation}
Also,
\begin{eqnarray*}
S_{n}(\theta) &=& - \Biggl[\frac{1}{n}\sum
_{k=1}^{n} \biggl(\frac{D^{1}_{\gamma}g_{\gamma
_{0}}}{g_{\gamma_{0}}}
(Y_{k} ) \biggr)^{2} \Biggr] \biggl[\sqrt {n}
\frac{q(\gamma_{1}-\gamma^{0})+p(\gamma_{2}-\gamma^{0})}{p+q} \biggr]^{2}
\\
&&{}+ \Biggl[\frac{1}{\sqrt{n}}\sum_{k=1}^{n}
\frac{D^{2}_{\gamma^{2}}g_{\gamma
_{0}}}{g_{\gamma_{0}}} (Y_{k} ) \Biggr] \biggl[\frac{q}{p+q}
\bigl(n^{1/4}\bigl(\gamma_{1}-\gamma^{0}\bigr)
\bigr)^{2}+\frac{p}{p+q} \bigl(n^{1/4}\bigl(
\gamma_{2}-\gamma^{0}\bigr) \bigr)^{2} \biggr]
\\
&&{}+ 2 \bigl(n^{1/4}\bigl(\gamma_{1}-\gamma^{0}
\bigr) \bigr)^{2} \Biggl[\frac{1}{\sqrt {n}}\sum
_{k=1}^{n}\bigl(D^{1}_{\gamma_{1}}p_{k}^{0}
\bigr)\frac{D^{1}_{\gamma
}g_{\gamma_{0}}}{g_{\gamma_{0}}} (Y_{k} ) \Biggr]
\\
&&{}-2 \bigl(n^{1/4}\bigl(\gamma_{2}-\gamma^{0}
\bigr) \bigr)^{2} \Biggl[\frac{1}{\sqrt {n}}\sum
_{k=1}^{n}\bigl(D^{1}_{\gamma_{2}}p_{k}^{0}
\bigr)\frac{D^{1}_{\gamma
}g_{\gamma_{0}}}{g_{\gamma_{0}}} (Y_{k} ) \Biggr]
\\
&&{}+2 \bigl(n^{1/4}\bigl(\gamma_{1}-\gamma^{0}
\bigr) \bigl(\gamma_{2}-\gamma^{0}\bigr) \bigr) \Biggl[
\frac{1}{\sqrt{n}}\sum_{k=1}^{n}
\bigl(D^{1}_{\gamma
_{2}}p_{k}^{0}-D^{1}_{\gamma_{1}}p_{k}^{0}
\bigr)\frac{D^{1}_{\gamma}g_{\gamma
_{0}}}{g_{\gamma_{0}}} (Y_{k} ) \Biggr].
\end{eqnarray*}
Using (\ref{recursion}) one gets that for all integer $k\geq2$
($D^{1}_{\gamma_{1}}p_{1}^{0}=0$ and $D^{1}_{\gamma_{2}}p_{1}^{0}=0$):
\[
D^{1}_{\gamma_{1}}p_{k}^{0}=
\frac{pq}{(p+q)^{2}}\sum_{l=1}^{k-1}(1-p-q)^{k-l}
\frac{D^{1}_{\gamma}g_{\gamma_{0}}}{g_{\gamma
_{0}}} (Y_{l} )
\]
and
\[
D^{1}_{\gamma_{2}}p_{k}^{0}=-D^{1}_{\gamma_{1}}p_{k}^{0}
\]
which leads to
\[
E_{\theta_{0}} \Biggl[ \Biggl( \frac{1}{\sqrt{n}}\sum
_{k=1}^{n}\bigl(D^{1}_{\gamma_{1}}p_{k}^{0}
\bigr)\frac{D^{1}_{\gamma}g_{\gamma
_{0}}}{g_{\gamma_{0}}} (Y_{k} ) \Biggr)^{2} \Biggr] \leq
\biggl( E_{\theta_{0}} \biggl(\frac{D^{1}_{\gamma}g_{\gamma_{0}}}{g_{\gamma
_{0}}} (Y_{1} )
\biggr)^{2} \biggr)^{2}
\]
and
\[
E_{\theta_{0}} \Biggl[ \Biggl( \frac{1}{\sqrt{n}}\sum
_{k=1}^{n}\bigl(D^{1}_{\gamma_{2}}p_{k}^{0}
\bigr)\frac{D^{1}_{\gamma}g_{\gamma
_{0}}}{g_{\gamma_{0}}} (Y_{k} ) \Biggr)^{2} \Biggr] \leq
\biggl( E_{\theta_{0}} \biggl(\frac{D^{1}_{\gamma}g_{\gamma_{0}}}{g_{\gamma
_{0}}} (Y_{1} )
\biggr)^{2} \biggr)^{2}.
\]
Thus, we obtain
%
%e36 #&#
\begin{equation}
\label{ordre2} \sup_{\theta\in U_{n}}\bigl\llvert S_{n}(
\theta) \bigr\rrvert =\mathrm{O}_{\Pb_{\theta
_{0}}} (1 ).
\end{equation}
For the order $3$ term, as soon as $\theta\in U_{n}$:
\begin{eqnarray*}
T_{n}(\theta) &=& - \Biggl[\sum_{k=1}^{n}
\biggl(\frac{D^{1}_{\gamma}g_{\gamma_{0}}}{g_{\gamma
_{0}}} (Y_{k} ) \biggr)^{3} \Biggr]
\biggl[\frac{q(\gamma
_{1}-\gamma^{0})+p(\gamma_{2}-\gamma^{0})}{p+q} \biggr]^{3}
\\
&&{}+ \Biggl[\sum_{k=1}^{n}
\frac{D^{3}_{\gamma^3}g_{\gamma_{0}}}{g_{\gamma
_{0}}} (Y_{k} ) \Biggr] \biggl[\frac{q}{p+q} \bigl(
\gamma_{1}-\gamma^{0} \bigr)^{3}+
\frac
{p}{p+q} \bigl(\gamma_{2}-\gamma^{0}
\bigr)^{3} \biggr]
\\
&&{}-3 \Biggl[\sum_{k=1}^{n}
\frac{D^{1}_{\gamma}g_{\gamma_{0}}}{g_{\gamma
_{0}}} (Y_{k} )\frac{D^{2}_{\gamma^{2}}g_{\gamma
_{0}}}{g_{\gamma_{0}}} (Y_{k} )
\Biggr] \biggl[\frac{q(\gamma
_{1}-\gamma^{0})+p(\gamma_{2}-\gamma^{0})}{p+q} \biggr]
\\
&&\hphantom{-}{}\times \biggl[\frac{q}{(p+q)^{2}} \bigl(\gamma_{1}-
\gamma^{0} \bigr)^{2}+\frac
{p}{(p+q)^{2}} \bigl(
\gamma_{2}-\gamma^{0} \bigr)^{2} \biggr]
\\
&&{}+\mathrm{O} \bigl(n^{-3/4} \bigr) \Biggl\{\sum
_{k=1}^{n}\bigl(D^{1}_{\gamma_{1}}p_{k}^{0}
\bigr) \biggl(\frac
{D^{1}_{\gamma}g_{\gamma_{0}}}{g_{\gamma_{0}}} (Y_{k} ) \biggr)^{2}
\\
&&\phantom{+\mathrm{O} \bigl(n^{-3/4} \bigr) \Biggl\{} {}+ \sum
_{k=1}^{n}\bigl(D^{1}_{\gamma_{1}}p_{k}^{0}
\bigr)\frac{D^{2}_{\gamma
^{2}}g_{\gamma_{0}}}{g_{\gamma_{0}}} (Y_{k} ) + \sum_{k=1}^{n}
\bigl(D^{2}_{(\gamma_{1})^{2}}p_{k}^{0}\bigr)
\frac{D^{1}_{\gamma
}g_{\gamma_{0}}}{g_{\gamma_{0}}} (Y_{k} )
\\
&&\phantom{+\mathrm{O} \bigl(n^{-3/4} \bigr) \Biggl\{} {}+\sum
_{k=1}^{n}\bigl(D^{2}_{(\gamma_{2})^{2}}p_{k}^{0}
\bigr)\frac{D^{1}_{\gamma
}g_{\gamma_{0}}}{g_{\gamma_{0}}} (Y_{k} ) + \sum_{k=1}^{n}
\bigl(D^{2}_{(\gamma_{1},\gamma_{2})}p_{k}^{0}\bigr)
\frac
{D^{1}_{\gamma}g_{\gamma_{0}}}{g_{\gamma_{0}}} (Y_{k} )  \Biggr\}
\end{eqnarray*}
so that using assumptions (\ref{hypocor2})
\[
\sup_{\theta\in U_{n}}\bigl\llvert T_{n}(\theta) \bigr\rrvert
=\mathrm{O}_{\Pb_{\theta_{0}}} \bigl(n^{-1/4} \bigr)+\mathrm{O}_{\Pb_{\theta_{0}}}
(1 )+ \mathrm{O} \bigl(n^{-1/4} \bigr)Z_{n}
\]
with
\begin{eqnarray*}
Z_{n}&=&\frac{1}{\sqrt{n}} \sum_{k=1}^{n}
\biggl\{ \biggl[ \biggl(\frac{D^{1}_{\gamma}g_{\gamma
_{0}}}{g_{\gamma_{0}}} (Y_{k} ) \biggr)^{2}+
\frac{D^{2}_{\gamma
^{2}}g_{\gamma_{0}}}{g_{\gamma_{0}}} (Y_{k} ) \biggr] D^{1}_{\gamma_{1}}p_{k}^{0}
\\
&&\phantom{\frac{1}{\sqrt{n}} \sum_{k=1}^{n}
\biggl\{} {}+\frac{D^{1}_{\gamma}g_{\gamma_{0}}}{g_{\gamma_{0}}} (Y_{k} ) \bigl[ D^{2}_{(\gamma_{1})^{2}}p_{k}^{0}+D^{2}_{(\gamma
_{2})^{2}}p_{k}^{0}+D^{2}_{(\gamma_{1},\gamma_{2})}p_{k}^{0}
\bigr]  \biggr\}.
\end{eqnarray*}
Now using (\ref{recursion}) one gets that for all integer $k\geq1$,
\begin{eqnarray*}
\frac{1}{1-p-q}D^{2}_{(\gamma_{1})^{2}}p_{k+1}^{0}&=&-2
\frac
{pq^{2}}{(p+q)^{3}} \biggl(\frac{D^{1}_{\gamma}g_{\gamma_{0}}}{g_{\gamma
_{0}}} (Y_{k} )
\biggr)^{2}
\\
&&{}+2\bigl(D^{1}_{\gamma_{1}}p_{k}^{0}
\bigr)\frac{D^{1}_{\gamma}g_{\gamma
_{0}}}{g_{\gamma_{0}}} (Y_{k} ) +\frac{pq}{(p+q)^{2}}\frac{D^{2}_{\gamma^{2}}g_{\gamma_{0}}}{g_{\gamma
_{0}}}
(Y_{k} ) +D^{2}_{(\gamma_{1})^{2}}p_{k}^{0},
\\
\frac{1}{1-p-q}D^{2}_{(\gamma_{2})^{2}}p_{k+1}^{0}&=&2
\frac
{p^{2}q}{(p+q)^{3}} \biggl(\frac{D^{1}_{\gamma}g_{\gamma_{0}}}{g_{\gamma
_{0}}} (Y_{k} )
\biggr)^{2}
\\
&&{}-2\bigl(D^{1}_{\gamma_{1}}p_{k}^{0}
\bigr)\frac{D^{1}_{\gamma}g_{\gamma
_{0}}}{g_{\gamma_{0}}} (Y_{k} ) -\frac{pq}{(p+q)^{2}}\frac{D^{2}_{\gamma^{2}}g_{\gamma_{0}}}{g_{\gamma
_{0}}}
(Y_{k} ) +D^{2}_{(\gamma_{2})^{2}}p_{k}^{0},
\\
\frac{1}{1-p-q}D^{2}_{(\gamma_{1},\gamma_{2})}p_{k+1}^{0}&=&2
\frac
{pq(q-p)}{(p+q)^{3}} \biggl(\frac{D^{1}_{\gamma}g_{\gamma_{0}}}{g_{\gamma
_{0}}} (Y_{k} )
\biggr)^{2}
\\
&&{}+2\bigl(D^{1}_{\gamma_{1}}p_{k}^{0}
\bigr)\frac{D^{1}_{\gamma}g_{\gamma
_{0}}}{g_{\gamma_{0}}} (Y_{k} ) +D^{2}_{(\gamma_{1},\gamma_{2})}p_{k}^{0},
\end{eqnarray*}
and using $D^{2}_{(\gamma_{1})^{2}}p_{1}^{0}=0, D^{2}_{(\gamma
_{2})^{2}}p_{1}^{0}=0, D^{2}_{(\gamma_{1},\gamma_{2})}p_{1}^{0}=0$ and
easy but tedious computations one gets that for some finite $C>0$,
\begin{eqnarray*}
&&E_{\theta_{0}} \bigl(Z_{n}^{2} \bigr)\leq C
E_{\theta_{0}} \biggl(\frac
{D^{1}_{\gamma}g_{\gamma_{0}}}{g_{\gamma_{0}}} (Y_{1} )
\biggr)^{2} \biggl[ E_{\theta_{0}} \biggl(\frac{D^{1}_{\gamma}g_{\gamma_{0}}}{g_{\gamma
_{0}}}
(Y_{1} ) \biggr)^{4}
\\
&&\phantom{E_{\theta_{0}} \bigl(Z_{n}^{2} \bigr)\leq C
E_{\theta_{0}} \biggl(\frac
{D^{1}_{\gamma}g_{\gamma_{0}}}{g_{\gamma_{0}}} (Y_{1} )
\biggr)^{2} \biggl[}{}+ E_{\theta_{0}} \biggl(\frac{D^{2}_{\gamma^{2}}g_{\gamma_{0}}}{g_{\gamma
_{0}}}
(Y_{1} ) \biggr)^{2} + \biggl(E_{\theta_{0}} \biggl(
\frac{D^{1}_{\gamma}g_{\gamma_{0}}}{g_{\gamma
_{0}}} (Y_{1} ) \biggr)^{2}
\biggr)^{2}   \biggr]
\end{eqnarray*}
so that we finally obtain
%
%e37 #&#
\begin{equation}
\label{ordre3} \sup_{\theta\in U_{n}}\bigl\llvert T_{n}(
\theta) \bigr\rrvert =\mathrm{O}_{\Pb_{\theta
_{0}}} (1 ).
\end{equation}
Let us finally study the fourth order remainder $R_{n}(\theta,t )$. We have
\[
\sup_{\theta\in U_{n}}\bigl\llvert R_{n}(\theta,t ) \bigr
\rrvert \leq\frac
{1}{n}\sum_{k=1}^{n}A_{k,n}B_{k,n},
\]
where, for big enough $n$, $A_{k,n}$ is a polynomial of degree at most
$4$ in $\sup_{\gamma'\in B_d(\gamma^{0}, \epsilon)}\| \frac
{D^{i}_{\gamma} g_{\gamma'}}{g_{\gamma'}}  (Y_{k} )\|$,
and $B_{k,n}$ is a sum of terms of form
%
%e38 #&#
\begin{equation}
\label{termB} \sup_{\theta\in U_{n}}\Biggl\llvert \prod
_{i=1}^{4}\prod_{j=0}^{i}
\bigl(D^{i}_{(\gamma_{1})^{j},(\gamma_{2})^{i-j}}p_{k}(\theta_{t})
\bigr)^{a_{i,j}}\Biggr\rrvert ,
\end{equation}
where the $a_{i,j}$ are non-negative integers such that $\sum_{i=1}^{4}\sum_{j=0}{i}a_{i,j} \leq4$.

To prove that
%
%e39 #&#
\begin{equation}
\label{ordre4} \sup_{\theta\in U_{n}}\bigl\llvert R_{n}(
\theta,t ) \bigr\rrvert =\mathrm{O}_{\Pb_{\theta
_{0}}} (1 )
\end{equation}
holds, it is enough to prove that
$E_{\theta_{0}}\llvert \sum_{k=1}^{n}A_{k,n}B_{k,n}\rrvert = \mathrm{O}(n)$. But
for each $k$, $p_{k}(\theta)$ and its derivatives depend on
$Y_{1},\ldots,Y_{k-1}$ only, so that $A_{k,n}$ and $B_{k,n}$ are
independent random variables, and
\begin{eqnarray*}
E_{\theta_{0}}\Biggl\llvert \sum_{k=1}^{n}A_{k,n}B_{k,n}
\Biggr\rrvert &\leq&\sum_{k=1}^{n}E_{\theta_{0}}
\llvert A_{k,n}\rrvert E_{\theta_{0}}\llvert B_{k,n}\rrvert
\\
&\leq&C \max_{i=1,2,3,4}E_{\theta_{0}} \biggl(\sup
_{\gamma'\in B_d(\gamma
^{0}, \epsilon)}\biggl\| \frac{D^{i}_{\gamma} g_{\gamma'}}{g_{\gamma'}} (Y_{1} )
\biggr\|^{4} \biggr) \sum_{k=1}^{n}E_{\theta_{0}}
\llvert B_{k,n}\rrvert
\end{eqnarray*}
for some finite $C>0$. Now,
using (\ref{recursion}) one gets that for all integer $k\geq 1$ and
for any $\theta$,
\begin{eqnarray*}
&&D^{1}_{\gamma_{1}}p_{k+1} (\theta )
\\
&&\quad= (1-p-q ) \biggl\{ \frac{p_{k}(\theta)(1-p_{k}(\theta))g_{\gamma_{2}}(Y_{k})D^{1}_{\gamma
}g_{\gamma_{1}}(Y_{k})+g_{\gamma_{1}}(Y_{k})g_{\gamma
_{2}}(Y_{k})D^{1}_{\gamma_{1}}p_{k} (\theta )} {
(p_{k}(\theta)g_{\gamma_{1}}(Y_{k})+(1-p_{k}(\theta))g_{\gamma
_{2}}(Y_{k}) )^{2}} \biggr\} ,
\\
&&D^{1}_{\gamma_{2}}p_{k+1} (\theta )
\\
&&\quad= (1-p-q ) \biggl\{ \frac{-p_{k}(\theta)(1-p_{k}(\theta))g_{\gamma_{1}}(Y_{k})D^{1}_{\gamma
}g_{\gamma_{2}}(Y_{k})+g_{\gamma_{1}}(Y_{k})g_{\gamma
_{2}}(Y_{k})D^{1}_{\gamma_{2}}p_{k} (\theta )} {
(p_{k}(\theta)g_{\gamma_{1}}(Y_{k})+(1-p_{k}(\theta))g_{\gamma
_{2}}(Y_{k}) )^{2}} \biggr\} .
\end{eqnarray*}
Notice that for any $\theta$, any $k\geq2$, $p_{k}(\theta) \in
(1-p,q)$ so that for any $\theta\in U_{n}$, any $k\geq2$, $p_{k}(\theta
) \in[\frac{1}{2}-\epsilon,\frac{1}{2}+\epsilon]$.
We obtain easily that for $i=1,2$, $k\geq2$,
\[
\sup_{\theta\in U_{n}} \bigl\llvert D^{1}_{\gamma_{i}}p_{k+1}
(\theta )\bigr\rrvert \leq \biggl(\frac{2\epsilon}{1-8\epsilon} \biggr) \biggl\{\sup
_{\gamma'\in
B_d(\gamma^{0}, \epsilon)}\biggl\| \frac{D^{1}_{\gamma} g_{\gamma'}}{g_{\gamma
'}} (Y_{k} )\biggr\|+ \sup
_{\theta\in U_{n}} \bigl\llvert D^{1}_{\gamma_{i}}p_{k}
(\theta )\bigr\rrvert \biggr\} .
\]
Using similar tricks, it is possible to get that there exists a finite
constant $C>0$ such that for any $i=1,2,3,4$, any $j=0,\ldots,i$, any
$k\geq2$,
\begin{eqnarray*}
&&\sup_{\theta\in U_{n}} \bigl\llvert D^{i}_{(\gamma_{1})^{j},(\gamma
_{2})^{i-j}}p_{k+1}(
\theta)\bigr\rrvert
\\
&&\quad\leq C \epsilon \Biggl\{\sup_{\gamma'\in B_d(\gamma^{0}, \epsilon)}\Biggl\| \sum
_{l=1}^{i}\frac
{D^{l}_{\gamma^{l}} g_{\gamma'}}{g_{\gamma'}} (Y_{k} )
\Biggr\|^{i+1-l}+ \sum_{l=1}^{i}\sum
_{m=0}^{l} \sup_{\theta\in U_{n}}
\bigl\llvert D^{l}_{(\gamma_{1})^{j},(\gamma
_{2})^{l-j}}p_{k}(\theta)\bigr\rrvert
^{i+1-l} \Biggr\}.
\end{eqnarray*}
By recursion, we obtain that there exists a finite $C>0$ such that any
term of form (\ref{termB}) has expectation uniformly bounded:
\begin{eqnarray*}
&&E_{\theta_{0}} \Biggl[\sup_{\theta\in U_{n}}\Biggl\llvert \prod
_{i=1}^{4}\prod
_{j=0}^{i} \bigl(D^{i}_{(\gamma_{1})^{j},(\gamma_{2})^{i-j}}p_{k}(
\theta _{t}) \bigr)^{a_{i,j}}\Biggr\rrvert \Biggr]
\\
&&\quad\leq C \max_{m=1,2,3,4} \max_{r=1,2,3,4}E_{\theta_{0}}
\biggl(\sup_{\gamma'\in
B_d(\gamma^{0}, \epsilon)}\biggl\| \frac{D^{m}_{\gamma} g_{\gamma'}}{g_{\gamma
'}} (Y_{1} )
\biggr\|^{r} \biggr)
\end{eqnarray*}
which concludes the proof of (\ref{ordre4}). Now,
using (\ref{mindeno1}), (\ref{Taylor1}), (\ref{ordre1}), (\ref
{ordre2}), (\ref{ordre3}) and (\ref{ordre4}), we get
\[
D_{n} \geq \mathrm{e}^{-\mathrm{O}_{\Pb_{\theta_{0}}} (1 )} \Pi_{2}
(U_{n} )
\]
so that \eqref{denom} holds with $S_n$ satisfying (\ref{Sn:D}).

%s4.4 #&#
\subsection{Proof of Theorem \texorpdfstring{\protect\ref{theo:gene}}{4}}
\label{subsec:gene}

The proof follows the same lines as in Ghosal and van~der Vaart
\cite{ghosal:vaart:2007}.
%Let $(\epsilon_n)_{n\geq1}$ be a sequence of positive real numbers as
%in assumption \textbf{C2} and \textbf{C3} or \textbf{C3bis}.
We write
\begin{eqnarray*}
&&\Pb^\Pi \biggl[ \|f_{l,\theta} - f_{l,\theta_0}
\|_1 \frac{\rho_\theta-
1}{2R_{\theta}+\rho_\theta- 1} \geq \epsilon_n \Big\vert
Y_{1:n} \biggr]
\\
&&\quad=\frac{ \int_{A_n \times{\mathcal X} } \mathrm{e}^{\ell_n(\theta,x) -\ell
_n(\theta_0,x_{0})}\Pi_{\Theta}(\mathrm{d}\theta) \pi_{\mathcal X} (\mathrm{d}x
) }{ \int_{\Theta\times{\mathcal X}} \mathrm{e}^{\ell_n(\theta,x) - \ell
_n(\theta_0,x_{0})}\Pi_{\Theta}(\mathrm{d}\theta) \pi_{\mathcal X} (\mathrm{d}x
) }
\\
&&\quad:= \frac{ N_n }{ D_n },
\end{eqnarray*}
where $A_n = \{ \theta\dvtx  \|f_{l,\theta} - f_{l,\theta_0}\|_1 \frac{\rho
_\theta- 1}{2R_{\theta}+\rho_\theta- 1} \geq \epsilon_n\}$.
A lower bound on $D_n $ is obtained in the following usual way. Set
$\Omega_n = \{(\theta, x) ; \ell_n(\theta,x) -\ell_n(\theta_0,x_{0})
\geq- n\tilde\epsilon_n^{2} \}$, which is a random subset of $\Theta
\times{\mathcal X}$ (depending on $Y_{1:n}$),% with $\tilde
\begin{eqnarray*}
D_n &\geq& \int_{S_n} \one_{\Omega_n}
\mathrm{e}^{\ell_n(\theta,x) -\ell_n(\theta
_0,x_{0})}\Pi_{\Theta}(\mathrm{d}\theta)
\pi_{\mathcal X} (\mathrm{d}x )
\\
&\geq& \mathrm{e}^{- n\tilde\epsilon_n^2 } \Pi(S_n \cap\Omega_n),
\end{eqnarray*}
therefore
%using assumption \textbf{C1} there exists $c_{1}>0$ such that
\begin{eqnarray*}
\Pb_{\theta_0} \bigl[ D_n <
\mathrm{e}^{-n\tilde\epsilon_n^2} \Pi(S_n)/2 \bigr] &\leq&
\Pb_{\theta_0} \bigl[ \Pi\bigl(S_n \cap\Omega_n^c
\bigr) \geq \Pi(S_n)/2 \bigr]
\\
&\leq& 2 \frac{ \int_{S_n} \Pb_{\theta_0} [\ell_n(\theta,x) - \ell
_n(\theta_0,x_{0}) \leq-n\tilde\epsilon_n^2  ] \Pi_{\Theta
}(\mathrm{d}\theta)\pi_{\mathcal X}(\mathrm{d}x)}{ \Pi(S_n ) }
\\
&=& \mathrm{o}(1)
\end{eqnarray*}
and
\[
\Pb^\Pi \biggl[ \|f_{l,\theta} - f_{l,\theta_0}
\|_1 \frac{\rho_\theta-
1}{2R_{\theta}+\rho_\theta- 1} \geq \epsilon_n \Big\vert
Y_{1:n} \biggr] = \mathrm{o}_{\Pb_{\theta_0}} (1 ) + \frac{N_{n}}{D_{n}}
\one_{2D_n \geq
\mathrm{e}^{-n\tilde\epsilon_n^2}\Pi(S_n)}.
\]
But
\begin{eqnarray*}
N_{n}&=&\int_{(A_n \cap{\mathcal F}_{n})\times{\mathcal X} } \mathrm{e}^{\ell
_n(\theta,x) -\ell_n(\theta_0,x_{0})}
\Pi_{\Theta}(\mathrm{d}\theta) \pi _{\mathcal X} (\mathrm{d}x )
\\
&&{}+\int_{(A_n \cap{\mathcal F}_{n}^{c})\times{\mathcal X} } \mathrm{e}^{\ell
_n(\theta,x) -\ell_n(\theta_0,x_{0})}
\Pi_{\Theta}(\mathrm{d}\theta) \pi _{\mathcal X} (\mathrm{d}x )
\end{eqnarray*}
and
\begin{eqnarray*}
&&E_{\theta_{0}} \biggl[\int_{(A_n \cap{\mathcal F}_{n}^{c})\times
{\mathcal X} } \mathrm{e}^{\ell_n(\theta,x) -\ell_n(\theta_0,x_{0})}
\Pi_{\Theta
}(\mathrm{d}\theta) \pi_{\mathcal X} (\mathrm{d}x ) \biggr]
\\
&&\quad=\mathrm{O} \bigl[\Pi_{\Theta} \bigl(A_n \cap{\mathcal
F}_{n}^{c} \bigr) \bigr]=\mathrm{o} \bigl(
\mathrm{e}^{-n\tilde\epsilon_n^2(C_n+1)} \bigr)
\end{eqnarray*}
by Fubini's theorem and assumption \textbf{C2} together with the fact
that $\ell_n(\theta_0)-\ell_n(\theta_0,x_{0})$ is uniformly upper
bounded. This implies
using assumption \textbf{C1} that
%
%e40 #&#
\begin{equation}
\label{intermediaire} \Pb^\Pi \biggl[ \|f_{l,\theta} -
f_{l,\theta_0}\|_1 \frac{\rho_\theta-
1}{2R_{\theta}+\rho_\theta- 1} \geq \epsilon_n
\Big\vert Y_{1:n} \biggr] =\mathrm{o}_{\Pb_{\theta_0}} (1 ) +
\frac{\tilde{N}_{n}}{D_{n}} \one _{2D_n \geq \mathrm{e}^{-n\tilde\epsilon_n^2}\Pi(S_n)},
\end{equation}
where $\tilde{N}_{n}=\int_{(A_n \cap{\mathcal F}_{n})\times{\mathcal
X} } \mathrm{e}^{\ell_n(\theta,x) -\ell_n(\theta_0,x_{0})}\Pi_{\Theta}(\mathrm{d}\theta)
\pi_{\mathcal X} (\mathrm{d}x ) $.
%
%(\ref{cond:KL:gene}).
Let now $(\theta_{j})_{j=1,\ldots, N}$, $N=N(\delta, \mathcal F_n,\allowbreak
d_{l}(\cdot,\cdot) )$, be the sequence of $\theta_j$'s in $\mathcal
F_n$ such for all $\theta\in\mathcal F_n$ there exists a $\theta_j$
with $d_{l}(\theta_j, \theta)\leq\delta$ with $\delta= \epsilon_n/12$.
Assume for simplicity's sake and without loss of generality that $n$ is
a multiple of the integer $l$, and define
\[
\phi_{j} = \one_{ \sum_{i=1}^{n/l} ( \one_{(Y_{li - l+1},\ldots,
Y_{li}) \in A_j } - \Pb_{\theta_0}((Y_{1},\ldots,Y_{l})\in A_j) )
> t_{j} },
\]
where
\[
A_j = \bigl\{ (y_{1},\ldots,y_{l}) \in
\mathcal Y^l \dvtx f_{l,\theta
_0}(y_{1},
\ldots,y_{l}) \leq f_{l,\theta_j}(y_{1},
\ldots,y_{l}) \bigr\}
\]
for some positive real number $t_{j}$ to be fixed later also.
Note that
\[
\Pb_{\theta_{j}}\bigl((Y_{1},\ldots,Y_{l})\in
A_j\bigr) -\Pb_{\theta
_0}\bigl((Y_{1},
\ldots,Y_{l})\in A_j\bigr) = \tfrac{ 1}{2 }
\|f_{l,\theta_j} - f_{l,\theta_0}\|_1.
\]
Define also
\[
\psi_n = \max_{1\leq j \leq N : \theta_{j}\in A_{n} }\phi_{j}.
\]
Then
%
%e41 #&#
\begin{equation}
\label{E0phin} E_{\theta_{0}} \biggl(\frac{\tilde{N}_{n}}{D_{n}}\psi_{n}
\biggr) \leq E_{\theta_{0}} \psi_{n} \leq N\bigl(\delta, \mathcal
F_n, d(\cdot,\cdot) \bigr)\max_{1\leq j \leq N :\theta_{j}\in A_{n}}
E_{\theta_{0}}\phi_{j}
\end{equation}
and
%
%e42 #&#
\begin{eqnarray}
\label{Ethetaphin} E_{\theta_{0}} \bigl(\tilde{N}_{n}(1-
\psi_{n}) \bigr) &=&\int_{\mathcal X} E_{\theta_{0},x_{0}}
\bigl(\tilde{N}_{n}(1-\psi_{n}) \bigr)\mu_{\theta
_{0}} (
\mathrm{d}x_{0} )
\nonumber
\\[-8pt]
\\[-8pt]
&=&\int_{(A_n \cap{\mathcal F}_{n})\times{\mathcal X}}E_{\theta,x} \bigl((1-\psi_{n})
\bigr)\Pi_{\Theta} (\mathrm{d}\theta )\pi_{\mathcal
X} (\mathrm{d}x ).
\nonumber
\end{eqnarray}
%
%so that:
%o\left(1\right)+ \left(\frac{n}{\delta}\right)^{M} \max_{1\leq j \leq
%N :\theta_{j}\in A_{n}}
%E_{\theta_{0}}\phi_{j} \\
%&+O\left[n^{D/2}e^{C_{n}}\int_{(A_n \cap{\mathcal F}_{n})\times{
%O\left[n^{D/2}\sup_{\theta\in A_{n}\cap{\mathcal F}_{n},x\in{
Now
\[
E_{\theta_0} [ \phi_j ] = \Pb_{\theta_0} \Biggl[ \sum
_{i=1}^{n/l} \bigl( \one_{(Y_{li - l+1},\ldots, Y_{li}) \in A_j } -
\Pb _{\theta_0}\bigl((Y_{1},\ldots,Y_{l})\in
A_j\bigr) \bigr) > t_{j} \Biggr]
\]
and
\begin{eqnarray*}
&&E_{\theta,x} (1-\phi_{j} )
\\
&&\quad = \Pb_{\theta,x}\Biggl [ \sum_{i=1}^{n/l}
\bigl( - \one_{(Y_{li - l+1},
\ldots,Y_{li}) \in A_j } + \Pb_{\theta,x}\bigl((Y_{li-l+1},
\ldots,Y_{li})\in A_j\bigr) \bigr)
\\
&&\phantom{\quad = \Pb_{\theta,x} \Biggl[} > -t_{j} +\sum
_{i=1}^{n/l} \bigl( \Pb_{\theta,x}
\bigl((Y_{li-l+1},\ldots ,Y_{li})\in A_j\bigr) -
\Pb_{\theta_0}\bigl((Y_{1},\ldots,Y_{l})\in
A_j\bigr) \bigr) \Biggr].
\end{eqnarray*}
Consider the sequence $(Z_{i})_{i\geq1}$ with for all $i\geq1$, $Z_i
= (X_{li-l+1}, \ldots,X_{li}, Y_{li-l+1}, \ldots,Y_{li})$, which is,
under $\Pb_\theta$, a Markov chain with transition kernel
%, which we denote by
$\bar{Q}_\theta$ given by
\begin{eqnarray*}
&&\bar{Q}_\theta\bigl(z,\mathrm{d}z'\bigr)
\\
&&\quad= g_{\theta}\bigl(y_1'|x_1'
\bigr) \cdots g_{\theta}\bigl(y_l'
|x_l'\bigr) Q_\theta\bigl(x_l,
\mathrm{d}x_1'\bigr) Q_\theta
\bigl(x_1',\mathrm{d}x_2'
\bigr)\cdots Q_\theta \bigl(x_{l-1}',
\mathrm{d}x_l'\bigr)\mu\bigl(\mathrm{d}y_1'
\bigr)\cdots\mu\bigl(\mathrm{d}y_l'\bigr).
\end{eqnarray*}
This kernel satisfies the same uniform ergodic property as $Q_\theta$,
with the same coefficients, that is condition \eqref{geo:ergo} holds
with the coefficients $R_{\theta}$ and $\rho_{\theta}$ with the
replacement of $Q_\theta$ by $\bar{Q}_\theta$, and we may
use Rio's \cite{rio:2000} exponential inequality (Corollary~1) with
uniform mixing coefficients (as defined in Rio \cite{rio:2000}) satisfying
$\phi(m) \leq R_\theta\rho_\theta^{-m}$. Indeed, by the Markov property,
\begin{eqnarray*}
\phi(m)&=&\sup_{A\in\sigma(Z_{1}), B\in\sigma(Z_{m+1})} \bigl(\Pb _{\theta}(B)-
\Pb_{\theta}(B \vert A) \bigr)
\\
&\leq& \sup_{z}\bigl | \Pb_{\theta}(Z_{m+1}\in B
) - \Pb_{\theta
}(Z_{m+1}\in B \vert Z_{1}=z )\bigr  |
\\
&\leq& R_\theta\rho_\theta^{-m}.
\end{eqnarray*}
We thus obtain that, for any positive real number $u$,
%
%e43 #&#
\begin{eqnarray}
\label{eq:rio1} &&\Pb_{\theta_0} \Biggl[ \sum_{i=1}^{n/l }
\bigl( \one_{(Y_{li -l+ 1},
\ldots,Y_{li}) \in A_j } - \Pb_{\theta_0}\bigl((Y_{1},
\ldots,Y_{l})\in A_j\bigr) \bigr) > u \Biggr]
\nonumber
\\[-8pt]
\\[-8pt]
&&\quad\leq\exp \biggl\{\frac{-2lu^{2} (\rho_{\theta_{0}}-1
)^{2}}{n (2R_{\theta_{0}}+\rho_{\theta_{0}}-1 )^{2}} \biggr\}
\nonumber
\end{eqnarray}
and
%
%e44 #&#
\begin{eqnarray}
\label{eq:rio2} &&\Pb_{\theta,x} \Biggl[ \sum_{i=1}^{n/l}
\bigl( - \one_{(Y_{li -l+ 1},
\ldots,Y_{li}) \in A_j } + \Pb_{\theta,x}\bigl((Y_{li-l+1},
\ldots,Y_{li})\in A_j\bigr) \bigr) > u \Biggr]
\nonumber
\\[-8pt]
\\[-8pt]
&&\quad\leq\exp \biggl\{\frac{-2lu^{2} (\rho_{\theta}-1 )^{2}}{n
(2R_{\theta}+\rho_{\theta}-1 )^{2}} \biggr\} .
\nonumber
\end{eqnarray}
Set now
\[
t_{j} = \frac{n\|f_{l,\theta_j} - f_{l,\theta_0}\|_1}{4l}.% \delta=
\]
Since for any $\theta$, $\frac{\rho_\theta- 1}{2R_{\theta}+\rho_\theta
- 1}\leq1$ and since consequently for $\theta_{j}\in A_{n}$, $\|
f_{l,\theta_j} - f_{l,\theta_0}\|_1 \geq\epsilon_{n}$, we first get,
using (\ref{eq:rio1}),
%
%e45 #&#
\begin{equation}
\label{eq:majo1} E_{\theta_0} [ \phi_j ]\leq\exp \biggl\{
\frac{-n \epsilon
_{n}^{2} (\rho_{\theta_{0}}-1 )^{2}}{8l (2R_{\theta
_{0}}+\rho_{\theta_{0}}-1 )^{2}} \biggr\}.
\end{equation}
Now, for any $\theta\in A_{n}$,
\begin{eqnarray*}
&& -t_{j} +\sum_{i=1}^{n/l}
\bigl( \Pb_{\theta,x}\bigl((Y_{li-l+1},\ldots ,Y_{li})\in
A_j\bigr) - \Pb_{\theta_0}\bigl((Y_{1},
\ldots,Y_{l})\in A_j\bigr) \bigr)
\\
&&\quad =-\frac{n\|f_{l,\theta_j} - f_{l,\theta_0}\|_1}{4l}+\frac{n}{l} \bigl\{ \Pb_{\theta_{j}}
\bigl((Y_{1},\ldots,Y_{l})\in A_j\bigr) -
\Pb_{\theta
_0}\bigl((Y_{1},\ldots,Y_{l})\in
A_j\bigr) \bigr\}
\\
&&\qquad{} + \frac{n}{l} \bigl\{\Pb_{\theta}\bigl((Y_{1},
\ldots,Y_{l})\in A_j\bigr) -\Pb_{\theta
_j}
\bigl((Y_{1},\ldots,Y_{l})\in A_j\bigr) \bigr
\}
\\
&&\qquad{}+\sum_{i=1}^{n/l} \bigl(
\Pb_{\theta,x}\bigl((Y_{li-l+1},\ldots,Y_{li})\in
A_j\bigr) - \Pb_{\theta}\bigl((Y_{1},
\ldots,Y_{l})\in A_j\bigr) \bigr)
\\
&&\quad \geq \frac{n\|f_{l,\theta_j} - f_{l,\theta_0}\|_1}{4l}-\frac{n\|f_{l,\theta
_j} - f_{l,\theta}\|_1}{l} -\sum
_{i=1}^{n/l}R_{\theta}\rho_{\theta}^{-i}
\\
&&\quad \geq\frac{n\|f_{l,\theta_j} - f_{l,\theta_0}\|_1}{4l}-\frac{n\|
f_{l,\theta_j} - f_{l,\theta}\|_1}{l} -\frac{R_{\theta} \rho_{\theta} }{\rho_{\theta}-1}
\\
&&\quad \geq\frac{n}{4l} \biggl(1- \frac{ 5 }{ 12} - \frac{ 4l }{ 12 n \epsilon
_n }
\biggr) \|f_{l,\theta} - f_{l,\theta_0}\|_1 \geq
\frac{n}{8l}\|f_{l,\theta} - f_{l,\theta_0}\|_1
\end{eqnarray*}
for large enough $n$, using the triangular inequality and the fact that
$\|f_{l,\theta_j} - f_{l,\theta}\|_1\leq\frac{ \epsilon_{n}}{12} \leq
\frac{\|f_{l,\theta} - f_{l,\theta_0}\|_1}{12} \frac{\rho_\theta-
1}{2R_{\theta}+\rho_\theta- 1}$
since $\theta\in A_{n}$ and $\frac{\rho_\theta- 1}{2R_{\theta}+\rho
_\theta- 1}\leq1$.
Then for $\theta\in A_{n}$ and large enough $n$,
%
%e46 #&#
\begin{equation}
\label{eq:majo2} E_{\theta,x} (1-\phi_{j} ) \leq\exp \biggl\{-
\frac{ n \epsilon_{n}^{2} }{ 32 l} \biggr\}.
\end{equation}
Combining \eqref{intermediaire}, with \eqref{E0phin}, \eqref{eq:majo1},
\eqref{Ethetaphin}, \eqref{eq:majo2} and using assumptions \textbf{C1}
and \textbf{C3} we finally obtain for large enough $n$
\begin{eqnarray*}
&& \Pb_{\theta_{0}} \biggl( \Pb^\Pi \biggl[ \|f_{l,\theta_j} -
f_{l,\theta
_0}\|_1 \frac{\rho_\theta- 1}{2R_{\theta}+\rho_\theta- 1} \geq \epsilon_n
\Big\vert Y_{1:n} \biggr] \biggr)
\\
&&\quad \leq \mathrm{o} (1 )+ \mathrm{O} \bigl(\mathrm{e}^{ n\tilde\epsilon_n^2 (1 + C_n) } \bigr)
\exp \biggl\{- \frac{ n \epsilon_{n}^{2} }{ 32 l} \biggr\}
\\
&&\qquad{}+ \exp \biggl\{\frac{-n \epsilon_{n}^{2} (\rho_{\theta_{0}}-1
)^{2}}{8l (2R_{\theta_{0}}+\rho_{\theta_{0}}-1 )^{2}} \biggr\} \exp \biggl\{
\frac{n \epsilon_n^2 ( \rho_{\theta_0}- 1)^2 }{ 16 l (
2R_{\theta_0} + \rho_{\theta_0}-1)^2 } \biggr\}
\\
&&\quad = \mathrm{o}(1).
\end{eqnarray*}
Assume now that assumption \textbf{C3bis} holds.
By writing $A_{n}\cap\mathcal F_n = \bigcup_{m\geq1} A_{n,m}(\epsilon
_{n})$ and using same reasoning, one gets, for some positive constant $c$:
\begin{eqnarray*}
&& \Pb_{\theta_{0}} \biggl( \Pb^\Pi \biggl[ \|f_{l,\theta_j} -
f_{l,\theta
_0}\|_1 \frac{\rho_\theta- 1}{2R_{\theta}+\rho_\theta- 1} \geq \epsilon_n
\Big\vert Y_{1:n} \biggr] \biggr)
\\
&&\quad = \mathrm{o} (1 )+ \mathrm{e}^{ n\tilde\epsilon_n^2 }\sum
_{m\geq1}\frac{ \Pi_{\Theta}
(A_{n,m}(\epsilon_{n}) ) }{ \Pi(S_n) } \exp \biggl\{- \frac{nm^{2}\epsilon_{n}^{2}}{ 32 l}
\biggr\}
\\
&&\qquad{} +\sum_{m\geq1} N \biggl(
\frac{m\epsilon_{n}}{12},A_{n,m}(\epsilon _{n}),d_{l}(
\cdot,\cdot) \biggr) \exp \biggl\{ -\frac{nm^{2} \epsilon_n^2 (
\rho_{\theta_0}- 1)^2 }{ 8 l ( 2R_{\theta_0} + \rho_{\theta_0}-1)^2
} \biggr\}
\\
&&\quad = \mathrm{o} (1 )
\end{eqnarray*}
and the second part of Theorem~\ref{theo:gene} is proved.
% \end{proof}

\section*{Acknowledgement}

The authors would like to thank the anonymous referees for providing
constructive comments that where very helpful to improve the
readability of the paper.
This work was partly supported by the 2010--2014 grant ANR Banhdits AAP
Blanc SIMI 1.

%suskaldyti doi

% imsref loaded by audrone.aklyte, 2014-01-08 08:57:29
%

\printhistory

\end{document}